\documentclass[reqno,11pt]{amsart}
\usepackage{geometry}
\usepackage{mathtools,amssymb,amsthm,mathrsfs,color,lineno,paralist,float}
\usepackage[colorlinks,
linkcolor=red,
anchorcolor=green,
citecolor=blue, 
]{hyperref}
\usepackage[T1]{fontenc}
\usepackage[utf8]{inputenc}

\numberwithin{equation}{section} 
\usepackage{calc}
\definecolor{bleu1}{RGB}{0,57,128}
\def\bleu1{\color{bleu1}}

\usepackage{etoolbox}
\patchcmd{\section}{\normalfont}{\normalfont \bleu1}{}{}
\patchcmd{\subsection}{\normalfont}{\normalfont \bleu1}{}{}
\patchcmd{\subsubsection}{\normalfont}{\normalfont \bleu1}{}{}

\newtheorem{THEalpha}{\bleu1 Theorem}

\newtheorem{The}{\bleu1 Theorem}[section]
\newtheorem{Cor}{\bleu1 Corollary}[section]
\newtheorem{Lem}{\bleu1 Lemma}[section]
\newtheorem{Pro}{\bleu1 Proposition}[section]
\theoremstyle{definition}
\newtheorem{defn}{\bleu1 Definition}[section]
\newtheorem*{Question}{Question}
\newtheorem{Problem}{Problem}
\theoremstyle{remark}
\newtheorem{Rem}{\bleu1 Remark}[section]

\newcommand{\TT}{\mathbb{T}}
\newcommand{\RR}{\mathbb{R}}
\newcommand{\ZZ}{\mathbb{Z}}
\newcommand{\NN}{\mathbb{N}}
\newcommand{\QQ}{\mathbb{Q}}
\newcommand{\CC}{\mathbb{C}}

\newcommand{\TA}{\mathcal{T}_{A,0}}
\newcommand{\TB}{\mathcal{T}_{B,0}}
\newcommand{\FF}{\mathbf{F}}
\newcommand{\GG}{\mathbf{G}}

\newcommand{\cT}{\mathcal{T}}
\newcommand{\cF}{\mathcal{F}}
\newcommand{\cD}{\mathcal{D}}
\newcommand{\cE}{\varepsilon}
\newcommand{\cU}{\delta}
\newcommand{\cL}{\mathcal{L}}
\newcommand{\cP}{\mathcal{P}}
\newcommand{\calE}{\mathcal{E}}
\newcommand{\cM}{\mathcal{H}}
\newcommand{\cV}{\mathcal{V}}

\newcommand{\bff}{\mathbf{f}}
\newcommand{\bfg}{\mathbf{g}}
\newcommand{\bfh}{\mathbf{h}}

\newcommand{\rS}{\mathrm{S}}
\newcommand{\rR}{\mathrm{R}}
\newcommand{\vep}{\varepsilon}
\newcommand{\expnx}{e^{i2\pi\langle n,x\rangle}}

\def\leq{\leqslant}
\def\geq{\geqslant}
\def\tilde{\widetilde}
\def\hat{\widehat}

\title[]{Rigidity properties for some isometric extensions of partially hyperbolic actions on the torus}
\author{Qinbo Chen}
\address[Qinbo Chen]{Department of Mathematics, Nanjing University, Nanjing 210093, China}
\email{qinbochen1990@gmail.com }
\author{Danijela Damjanovi\'{c}}
\address[Danijela Damjanovi\'{c}]{Department of Mathematics, Kungliga Tekniska H\"{o}gskolan, Lindstedtsv\"{a}gen 25, SE-100 44  Stockholm, Sweden}
\email{ddam@kth.se}


\begin{document} 
\begin{abstract}
This paper studies local rigidity for some isometric toral extensions of partially hyperbolic $\mathbb{Z}^k$ ($k\geqslant 2$) actions on the torus. We prove a $C^\infty$ local rigidity result for such actions, provided that the smooth perturbations of the actions satisfy the intersection property. We also give a local rigidity result within a class of volume preserving actions. Our method mainly uses a generalization of the KAM iterative scheme.
\end{abstract}
\maketitle


\section{Introduction}
Let $A$ be an automorphism of the  torus  $\TT^d=\RR^d/\ZZ^d$. For a $C^\infty$ function $\tau(x): \TT^d\longrightarrow\RR^s$ with the integer $s\geq 1$, it defines an isometric toral extension of $A$, which is a map $\cT_{A,\tau}:\TT^d\times\TT^s\longrightarrow \TT^d\times\TT^s$ of the form
 \begin{equation}\label{gen_isoextension}
	\cT_{A,\tau}(x,y)=(Ax,y+\tau(x) \textup{~mod~}\ZZ^{s})
\end{equation}
where  $\TT^s=\RR^s/\ZZ^s$. We can think of $\TT^d\times\TT^s$ as a (trivial) bundle over the base space $\TT^d$, so that $\cT_{A,\tau}$ is a skew product with an automorphism $A$ on the base $\TT^d$ and a translation on each fiber $\{x\}\times\TT^s$ with the translation vector $\tau(x)$.  If $A$ is ergodic, then such isometric extensions provide a special class of volume preserving partially hyperbolic systems. They have been extensively studied in the literature, especially in the case when $A$ is hyperbolic.

Our paper treats abelian group actions generated by multiple diffeomorphisms of the form \eqref{gen_isoextension}, and studies the rigidity properties of such actions under perturbations. In general, a smooth $\ZZ^k$ action $\rho$ on a compact nilmanifold (including the torus)  $M$ is given by a group morphism  $\rho:\mathbf{n}\mapsto  \rho(\mathbf{n})$ from  $\ZZ^k$  into the group $\textup{Diff}^\infty(M)$ of $C^\infty$ diffeomorphisms of $M$. The classification of smooth actions of higher rank on compact manifolds is one of the central problems in smooth dynamics. It originated from  the Zimmer program of studying actions of higher rank groups and lattices \cite{Zimmer1987}.  
 
 The action considered in this paper acts by automorphisms on the base $\TT^d$ and acts isometrically on the fiber $\TT^s$. We start by considering a class of $\ZZ^2$ actions $\alpha=\langle\cT_{A_1,\tau_1}, \cT_{A_2,\tau_2}\rangle$, that is $\alpha(\mathbf{n})=\cT^{n_1}_{A_1,\tau_1}\circ \cT^{n_2}_{A_2,\tau_2}$ for all $\mathbf{n}=(n_1, n_2)\in \ZZ^2$. We are motivated by an attempt to understand the smooth actions close to $\alpha$ in terms of their dynamics and  geometry. Moreover, we assume that the $\ZZ^2$ action $\langle A_1, A_2\rangle$ on the base $\TT^d$ is higher rank, which is closely tied to certain ergodic properties (see Remark \ref{Rem_ergod}). As a consequence, one can show that $ \cT_{A_i,\tau_i}$, $i=1,2$, are simultaneously $C^\infty$-conjugate to $ \cT_{A_i,[\tau_i]}$, see Proposition \ref{Pro_conj_ave}, where the constant vector
    \[[\tau_i]\overset{\Delta}=\int_{\TT^d}\tau_i(x)\, dx\in \RR^s\] 
    denotes the average of $\tau_i$ over $\TT^d=\RR^d/\ZZ^d$, $i=1,2$.

People are interested in the smooth rigidity problem of the above actions. By the discussion above, it is related to the properties of the averages $[\tau_1]$ and $[\tau_2]$. We first recall some prior works. 

 $\bullet$ If $s=0$ (i.e. no any extensions), then $\alpha$ becomes $\langle A_1, A_2\rangle$, and the local rigidity has been established by Katok and the second author. More precisely, 
\begin{The}\label{thmDamjKatok10}\cite{Damjanovic_Katok10}
If the $\ZZ^2$ action $\alpha=\langle A_1, A_2\rangle$ is higher rank, then there exists an integer $l=l(\alpha)>0$ such that: any smooth action $\tilde\alpha:\ZZ^2\to \textup{Diff}^\infty(\TT^d)$ which is  sufficiently close to $\alpha$ in the $C^l$ topology is $C^\infty$-conjugate to $\alpha$. 
\end{The}

It still holds for any higher rank $\ZZ^k$, $k\geq 2$ actions by toral automorphisms, cf. \cite{Damjanovic_Katok10}.

$\bullet$ If $s\geq 1$, the action may enjoy a local rigidity subject to constraints that some invariants are  preserved. For $[\tau_1]$ and $[\tau_2]$ satisfying certain Diophantine condition, the following form of local rigidity holds.
\begin{The}\cite{Damjanovic_Fayad}
Consider the $\ZZ^2$ action $\alpha=\langle \cT_{A_1, \tau_1},\cT_{A_2,\tau_2}\rangle$, where the averages $[\tau_1]$ and $[\tau_2]$  satisfy the simultaneous Diophantine condition and the action $\langle A_1, A_2\rangle$ on the base $\TT^d$ is higher rank. Then, there exists an integer $l=l(\alpha)>0$ such that: for any smooth $\ZZ^2$ action $\tilde{\alpha}=\langle F_1, F_2 \rangle$	that is sufficiently close to $\alpha$ in the $C^l$ topology, 
    $\tilde\alpha$ can be $C^\infty$-conjugate to $\alpha$ provided that each $F_i$, $i=1,2$, preserves an invariant probability measure $\mu_i$ whose translation vector along the fiber direction is equal to $[\tau_i]$. 
\end{The}
The above assumption is inspired by Moser's local rigidity result \cite{Moser90_commuting} for commuting circle maps (see also \cite{Fayad_Khanin} for a global result).

Nevertheless, very little is known for the rational case, i.e., $[\tau_i]\in \QQ^s$, $i=1,2$. In this situation, $\cT_{A_i,\tau_i}$, $i=1,2$ are \textbf{non-ergodic} on the total space $\TT^d\times\TT^s$, thus in order to study the rigidity aspect of such actions, some  constraints are needed to be imposed on the class of perturbations. More precisely, one may ask the following question.   
   
\begin{Question}
Let $[\tau_1]\in \QQ^s$ and $[\tau_2]\in \QQ^s$, and consider  a smooth $\ZZ^2$ action $\tilde\alpha$  that is close to $\alpha=\langle \cT_{A_1,\tau_1}, \cT_{A_2,\tau_2}\rangle$  in the $C^r$ topology with $r$ being suitably large. Under which conditions can we show that the perturbed action $\tilde\alpha$ is $C^\infty$-conjugate to the original action $\alpha$?
\end{Question}

In this paper, we solve this problem by only assuming that the perturbed actions  satisfy certain topological assumption (see subsection \ref{subsecmainres}). This assumption is not only sufficient but also necessary. 
   
We stress that the actions considered here are not necessarily Anosov on the base $\TT^d$, and therefore, they are not the so-called fibered partially hyperbolic systems discussed in \cite{Damjanovic_Wilkinson_Xu2021_Duke}. In the case of Anosov base, there are many geometric tools that can be used towards classifying perturbations, see e.g. \cite{Damjanovic_Wilkinson_Xu2021_Duke} and references therein.

We also remark that for groups with more structure than $\ZZ^k$, the perturbations are better understood. In \cite{FM_2009}, Fisher and Margulis established local rigidity in full generality for quasi-affine actions by higher rank lattices in semisimple Lie groups. Prior to \cite{FM_2009}, the question about local rigidity of product actions of property (T) groups  has been addressed in  \cite{Nitica_Torok_1995,Nitica_Torok_2001,Torok_2003} where they considered   higher-rank lattice actions of the form $\rho\times id_{\TT^1}$ with the subaction $\rho$ having certain hyperbolic structure.  However, the situations and methods in these works are very different from ours and depend on the acting group having Kazhdan's property (T).

\subsection{Background on rigidity for smooth actions}
To better explain the background and motivation of our result, we give a brief introduction to the rigidity problem of smooth actions mainly from  the viewpoint of dynamical systems. The interested reader can also refer to \cite{Fishier_2007Local} for a survey of the local rigidity problem for general group actions.  Let $M$ be a compact manifold. We refer to a homomorphism $\rho: \ZZ^k\to \textup{Diff}^\infty(M)$ as an action since it can be thought of as $C^\infty$ action $\rho: \ZZ^k\times M\to M$.
Briefly, we say a $\ZZ^k$ action $\rho$ is  $C^\infty$-locally rigid if for any sufficiently small perturbations $\tilde\rho$, there is a $C^\infty$ conjugacy $h$ such that $h\circ \tilde\rho(\mathbf{n})\circ h^{-1}= \rho(\mathbf{n})$, for all $\mathbf{n}\in\ZZ^k$. Two $\ZZ^k$ actions $\rho$ and $\tilde\rho$ are said to be $C^r$-close if the diffeomorphisms $\rho(\mathbf{e}_i)$ and $\tilde\rho(\mathbf{e}_i)$ are close in the $C^r$ topology for all $i=1,\cdots, k$, where $\mathbf{e}_1,\cdots,\mathbf{e}_k$ are the generators of $\ZZ^k$.

The dynamical motivation for investigating the rigidity started with the study of structural stability. In contrast to the  structural stability which preserves only the topological orbit structure,  the $C^\infty$-local rigidity preserves all differentiable orbit  structure. By a classical result of Franks and Manning, any  Anosov diffeomorphism on the torus is topologically conjugate to an affine Anosov automorphism. Nevertheless,  this topological conjugacy, generically, cannot be improved to $C^1$. For example, one can easily perturb an Anosov diffeomorphism $f$ around a periodic point to a new Anosov diffeomorphism $\tilde f$ which changes the eigenvalues of its differential at this periodic point, and thus $\tilde f$  cannot be $C^1$ conjugate to $f$.

In contrast to the rank-one (i.e. $\ZZ^1$ action) situation,  higher rank abelian Anosov actions exhibit much more rigidity. Some special cases are studied in \cite{Katok_Lewis_1991,Hur_1992,Katok_Lewis_Zimmer_1996}. Katok and Spatzier first established the local rigidity for  higher-rank algebraic Anosov actions with semisimple linear parts \cite{KS_1997}, and later extended to some non-semisimple actions by Einsiedler and Fisher \cite{Einsiedler_Fisher2007}. These results motivate the Katok-Spatzier global rigidity conjecture: \textit{all irreducible higher rank abelian smooth Anosov actions on a compact manifold are smoothly conjugate to algebraic actions}. It is concerned with classification of higher-rank Anosov smooth actions. Over the last two decades, significant progress has been made  towards this conjecture, we only list a few \cite{Kal_Sa_2006,Kalinin_Spatzier2007,RodriguezHertz_globalrigidity_2007,FKS_2013,RW_2014} and  see references therein. In particular,  Rodriguez Hertz and Wang \cite{RW_2014}  obtained the optimal global rigidity  result on nilmanifolds and tori for $\ZZ^k$ Anosov actions without rank-one factors. This extended earlier work \cite{FKS_2013} which required every Weyl chamber contains an Anosov element. 

However, smooth classification of partially hyperbolic actions is much more complicated. Even local rigidity results are scarce. The major difficulty lies in the appearance of center foliations. The first breakthrough is the work \cite{Damjanovic_Katok10} on local rigidity of certain partially hyperbolic affine actions of abelian groups on tori using a KAM approach. A few more recent developments and rigidity results on partially hyperbolic actions can be found in \cite{Damjanovic_Katok2011,Damjanovic_Fayad,Vin_Wang_2019, Damjanovic_Wilkinson_Xu2021_Duke}, \textit{etc}.

\subsection{The main results} \label{subsecmainres}

The actions considered in this paper are a special class of partially hyperbolic actions: (I) the action has a partially hyperbolic part and a non-hyperbolic, isometric extension part; (II) all elements of the action are non-ergodic on the total space.

To state our main results we need some basic definitions. The first one is the higher rank condition which is a common condition used in the rigidity problem of group actions. 
\begin{defn}
We say that a $\ZZ^k$, $k\geq 2$ action has \emph{a rank-one factor}  if it factors to a $\ZZ^k$ action which is (up to a finite index subgroup of $\ZZ^k$) generated by a single diffeomorphism. Moreover, an action is said to be \emph{higher rank} if it has no rank-one factors.
\end{defn}

\begin{Rem}[Ergodicity]\label{Rem_ergod}
From a dynamical viewpoint, a $\ZZ^k$ action by toral automorphisms is higher rank is equivalent to saying that it contains a subgroup $L$ isomorphic to $\ZZ^2$ such that every element in $L$, except for identity, is ergodic. See \cite{Starkov1999}.  In particular,  for a $\ZZ^2$ action, this condition is equivalent to saying that all non-trivial elements of the action are ergodic, see Lemma \ref{lem_equihighrank} for an explanation. Consequently, these ergodic automorphisms are partially hyperbolic. 
\end{Rem}

For our purpose we need  the following notion.
\begin{defn}
A map $F(x,y): \TT^d\times\TT^s\to  \TT^d\times\TT^s$  is said to satisfy \emph{the intersection property} if the following condition {\bf(IP)} holds:
	\begin{enumerate}
    \renewcommand{\labelenumi}{\theenumi}
    \renewcommand{\theenumi}{\bf(IP)}
    \makeatletter
    \makeatother
	\item  \label{condIP}  for any $d$-dimensional subtorus $\Gamma$ which is diffeomorphic and $C^1$-close\footnote{It means that $\Gamma$ is a $d$-dimensional submanifold of $\TT^d\times\TT^s$, and it is of the form $\{(x,y): x\in \TT^d,\quad y=y_0+\psi(x) ~\textup{mod}~\ZZ^s\}$ where $\psi\in C^1(\TT^d,\RR^s)$ and $\|\psi\|_{C^1}\leq \delta$ for a priori fixed number $\delta>0$.} to $\TT^d\times\{y_0\}$, with $y_0\in \TT^s$,  one has $F(\Gamma)\cap \Gamma\neq \emptyset.$ 
    \end{enumerate}
\end{defn}

For instance, the map $\cT_{A,0}=A\times id_{\TT^s}$ satisfies condition \ref{condIP}. This can be readily verified by using the fact that $x=0$ is always a fixed point of the automorphism $A:\TT^d\to\TT^d$.  

Historically, the intersection property was once used by Moser (see  \cite{Siegel_Moser1971} or   R\"ussmann's works), as an alternative for the area-preserving condition, to prove the existence of  invariant circles for the twist maps of a cylinder (Moser's twist map theorem).  Subsequently,  several high-dimensional versions of the intersection property were introduced to study the existence of invariant tori for certain non-symplectic maps in high-dimensional spaces. In the present paper condition \ref{condIP} also plays an important role in the proof.

Now we can state the first main result. 

\begin{THEalpha}\label{MainThm_0}
Consider the smooth $\ZZ^2$ action $\alpha=\langle\cT_{A_1,\tau_1}, \cT_{A_2,\tau_2}\rangle$ on $\TT^d\times\TT^s$, where the averages $[\tau_1]$ and $[\tau_2]$ are rational and the action $\langle A_1, A_2\rangle$ on the base $\TT^d$ is higher rank. Then, there exist  $\vep=\vep(\alpha)>0$ and $\mu=\mu(\alpha)>0$ such that: for any smooth $\ZZ^2$ action  $\tilde{\alpha}=\langle\cF_1, \cF_2\rangle$, if 
\begin{enumerate}[(i)]
	\item $\textup{dist}_{C^\mu}(\tilde{\alpha}, \alpha)<\vep$.	
	\item the finite set $\{(i,j)\in \ZZ^2: ~|i|\leq M_0, |j|\leq M_0\}$ contains two linearly independent 
	elements $\mathbf{m}, \mathbf{n}$ such that $\tilde{\alpha}(\mathbf{m})$ and $\tilde{\alpha}(\mathbf{n})$ satisfy condition \ref{condIP}, where  $M_0\geq 1$ is the minimal positive integer $\lambda$ such that  $\lambda\,[\tau_1]\in \ZZ^s$ and $\lambda\,[\tau_2]\in \ZZ^s$. 
\end{enumerate}
then the action $\tilde{\alpha}=\langle \cF_1, \cF_2\rangle$ is $C^\infty$-conjugate to $\alpha=\langle\cT_{A_1,\tau_1}, \cT_{A_2,\tau_2}\rangle$.
\end{THEalpha}

\begin{Rem}
By assumption we see that the original action $\alpha$ acts ergodically on the base $\TT^d$ and acts isometrically on the fiber $\TT^s$. But $\alpha$ does not act ergodically on the total space $\TT^d\times\TT^s$ since $[\tau_1], [\tau_2]$ are rational. For the generators of the perturbed action $\tilde\alpha$, $\cF_l(x,y)\in \textup{Diff}^\infty(\TT^d\times\TT^s)$, $l=1,2$, can be written in the form $\cF_l=\cT_{A_l, \tau_l}+f_l$ where the perturbation term $f_l(x,y)=(f_{l,1}(x,y), f_{l,2}(x,y))$ with $f_{l,1}\in C^\infty(\TT^{d+s}, \RR^d)$ and $f_{l,2}\in C^\infty(\TT^{d+s}, \RR^s)$. The above mentioned $C^\mu$ distance between two actions $\tilde\alpha$ and $\alpha$ is defined by using the generators:
\[\textup{dist}_{C^\mu}(\tilde{\alpha}, \alpha):=\max\limits_{l=1,2}\|\cF_l-\cT_{A_l,\tau_l}\|_{C^\mu(\TT^d\times\TT^s)}.\]	
\end{Rem}

\begin{Rem}
We will see from the proof of Theorem \ref{MainThm_0} that there is a near-identity conjugacy $U\in \textup{Diff}^\infty(\TT^d\times\TT^s)$ such that $U\circ \tilde\alpha(\mathbf{k})\circ U^{-1}= \alpha(\mathbf{k})$, for all $\mathbf{k}\in \ZZ^2$. As will be shown in subsection \ref{subsecproofA}, the construction of the conjugacy $U$ mainly consists of two  parts: one is produced by the KAM scheme, and the other is obtained by solving a cohomology equation over periodic diffeomorphisms.   
\end{Rem}

\begin{Rem}
Condition \ref{condIP} cannot be removed, otherwise, the above result may fail. Although the local rigidity of the higher rank action $\langle A, B\rangle$ on $\TT^d$ holds (see Theorem \ref{thmDamjKatok10}), it can not be applied directly to prove the situation considered by the above theorem. This is because the generating elements $\cF_1(x,y), \cF_2(x,y)$ of the perturbed action $\tilde\alpha$ 	depend on both the base variable $x$ and the fiber variable $y$, and for each fixed $y$ the restriction map $\pi_1\circ\cF_1(\cdot,y), \pi_1\circ\cF_2(\cdot, y):$ $\TT^d\to \TT^d$ of $\cF_1$ and $\cF_2$ on the base $\TT^d$, generically, do not commute. Here,  $\pi_1:(x,y)\mapsto x$ is the projection.  
\end{Rem}

\subsubsection{Volume preserving actions}
By a volume preserving diffeomorphism of $M$, we mean a diffeomorphism that preserves a volume form on $M$. For the actions considered here, if \textbf{the fiber dimension $s=1$} and the perturbations are within the class of volume preserving actions, we have the following result.

\begin{THEalpha}\label{MainThm_1}
Consider the smooth $\ZZ^2$ action $\alpha=\langle\cT_{A_1,\tau_1}, \cT_{A_2,\tau_2}\rangle$ on $\TT^d\times\TT^1$, where the averages $[\tau_1]$ and $[\tau_2]$ are rational and the action  $\langle A_1, A_2\rangle$ on the base is higher rank. Then, there exist  $\vep=\vep(\alpha)>0$ and   $\mu=\mu(\alpha)>0$ such that: for any smooth $\ZZ^2$ action  $\tilde{\alpha}=\langle \cF_1, \cF_2\rangle$ with $\cF_l$, $l=1,2$, preserving a volume form on $\TT^{d+1}$, if
   \begin{enumerate}[(i)]
	\item $\textup{dist}_{C^\mu}(\tilde{\alpha}, \alpha)<\vep$.	
	\item for $l=1,2$ let $q_l$ denote the minimal positive integer $\lambda$ such that  $\lambda\,[\tau_l]\in \ZZ$, and assume that $\cF_l^{q_l}=\cF_l\circ\cdots\circ\cF_l$ admits an invariant $d$-dimensional torus homotopic to $\TT^d\times\{0\}$.
   \end{enumerate}
then the action $\tilde{\alpha}=\langle \cF_1, \cF_2\rangle$ is $C^\infty$-conjugate to $\alpha=\langle\cT_{A_1,\tau_1}, \cT_{A_2,\tau_2}\rangle$.
\end{THEalpha}

We also give a corresponding result for $\ZZ^k$, $k\geq 2$ actions. 
\begin{THEalpha}\label{MainThm_2}
Let $\rho=\rho_0\times id_{\TT^1}$ be a $\ZZ^k$, $k\geq 2$ action on $\TT^d\times\TT^1$, where $\rho_0$ is a $\ZZ^k$ higher rank action by automorphisms on $\TT^d$. Then, there exist  $\vep=\vep(\rho)>0$ and   $\mu=\mu(\rho)>0$ such that: for any smooth $\ZZ^k$ action  $\tilde{\rho}$ which preserves a volume form on $\TT^{d+1}$, if
   \begin{enumerate}[(i)]
	\item $\textup{dist}_{C^\mu}(\tilde{\rho}, \rho)<\vep$.	
	\item $\tilde{\rho}$ admits a common invariant $d$-dimensional torus homotopic to $\TT^d\times\{0\}$.
   \end{enumerate}
then  $\tilde{\rho}$ is $C^\infty$-conjugate to $\rho$.	
\end{THEalpha}
\begin{Rem}
In condition $(ii)$, $\tilde{\rho}$ admits a common invariant $d$-dimensional torus means that there is a $d$-dimensional torus which is invariant under $\tilde\rho({\mathbf{e}_i})$, for all generators $\mathbf{e}_1,\cdots, \mathbf{e}_k\in \ZZ^k$.
\end{Rem}

\subsubsection{The method.} Different from almost all existing local/global results for isometric extensions where the base maps are Anosov, the situation considered here is only \textit{partially hyperbolic} on the base $\TT^d$, so the geometric methods in previous works on isometric extensions of Anosov systems or isometric extensions of partially hyperbolic, accessible systems are not applicable in our situation. The main approach we will use is a generalization of the KAM (Kolmogorov-Arnold-Moser) method. In fact, finding a conjugacy between a perturbed action and the original one is a problem of inverting a nonlinear operator. Our strategy is to apply linearization and successive iterations to produce a solution to the nonlinear problem. It is independent of methods which use hyperbolic dynamics on the base, so our arguments are analytic rather than geometric in nature.

\subsection{Strategy of the proof}
The higher rank condition for the action $\langle A_1, A_2\rangle$ on the base $\TT^d$ plays a key role throughout  our proofs. One immediate consequence is that  $ \cT_{A_i,\tau_i}$, $i=1,2$, are simultaneously $C^\infty$-conjugate to $\cT_{A_i,[\tau_i]}$. The philosophy behind this fact is that the higher rank condition on the base implies that $\tau_i(x)$, $i=1,2$ is a coboundary with respect to the action $\langle A_1, A_2\rangle$, see Section \ref{Section_actionofproducttype}. Since  $\cT_{A_i,[\tau_i]}=A_i\times R_{[\tau_i]}$, if $q_i$ is the period of the rational vector $[\tau_i]$, then the $q_i$-fold composition becomes $\cT_{A_i, [\tau_i]}^{q_i}=A_i^{q_i}\times id_{\TT^s}$. Consequently, this reduces Theorem \ref{MainThm_0} to studying the local rigidity of a subaction generated by two automorphisms of the form $\cT_{A,0}=A\times id_{\TT^s}$ and $\cT_{B,0}=B\times id_{\TT^s}$, where $A, B$ still satisfy the higher rank condition. See Theorem \ref{Element_Thm1}.

The proof of Theorem \ref{Element_Thm1} is included in Section \ref{Section_conjugequation}--Section \ref{Section_KAMscheme}.  We adopt the KAM methodology, as in \cite{Damjanovic_Katok10}. The basic philosophy of this method is that one can reduce a nonlinear problem to a linear one and solve approximately the linearized equations, then by iterating this process, the limit of successive iterations finally produces a solution to the nonlinear problem. The proof includes obtaining tame solutions for the cohomological equations and constructing tame splitting as well. 

More precisely,  the linear problem consists of solving approximately (with error quadratically small with respect to the error of the perturbation) two kinds of cohomological equations over a $\ZZ^2$ action by non-ergodic partially hyperbolic automorphisms: the twisted case
\begin{equation}\label{sec1twis}
\begin{aligned}
	u_1(Ax, y)-Au_1(x,y)=\bff_1(x,y),\qquad u_1(Bx,y)-Bu_1(x,y)=\bfg_1(x,y),
\end{aligned}	
\end{equation}
and the untwisted case
\begin{equation}\label{sec1untwis}
  u_2 (Ax,y)-u_2(x,y)=\bff_2(x,y),\qquad u_2(Bx,y)-u_2(x,y)=\bfg_2(x,y).\
\end{equation}
If $\cL_1(\bff_1, \bfg_1)=0$ and $\cL_2(\bff_2, \bfg_2)=0$ (operators $\cL_1$ and $\cL_2$ are defined in \eqref{fas1}-\eqref{fas2}) are satisfied, then due to the higher rank condition it is feasible to solve equations \eqref{sec1twis}-\eqref{sec1untwis} exactly and estimate the solution with a fixed loss of regularity. For a general perturbation $\cL_1(\bff_1, \bfg_1)\ne 0$ and $\cL_2(\bff_2, \bfg_2)\ne 0$, one needs (as in \cite{Damjanovic_Katok10}) to approximate data given by the perturbation by data which satisfies $\cL_1(\bff_1, \bfg_1)=0$ and $\cL_2(\bff_2, \bfg_2)=0$. To this aim, we use the concrete constructions from \cite{Damjanovic_Katok10} but we have to combine them with the idea of smooth dependence on parameters to give tame estimates (for the solutions) in the fiber direction as well as the base direction. See Section \ref{Section_conjugequation} and Section \ref{Section_Smoothdependparamt} for more details.

When it comes to verifying the convergence of the KAM scheme, one needs to handle the following problems:\\
(I)\, $\cL_1(\bff_1, \bfg_1)\neq 0$ and $\cL_2(\bff_2, \bfg_2)\neq 0$ in general;\\
(II)\, the fixed loss of regularity; \\
(III)\, the estimate of the averaged terms $[\bff_2](y):=\int_{\TT^d}\bff_2(x,y)\, dx$ and $[\bfg_2](y):=\int_{\TT^d}\bfg_2(x,y)\,dx$.\\ 
To tackle (I), the basic idea is to solve equations \eqref{sec1twis}--\eqref{sec1untwis} up to quadratic errors. This requires to split $\bff_i, \bfg_i$, $i=1,2$  into $\bff_i=\cP(\bff_i)+\calE(\bff_i)$ and $\bfg_i=\cP(\bfg_i)+\calE(\bfg_i)$ in a tame way, such that $\cL_i(\cP(\bff_i), \cP(\bfg_i))=0$ and the remainder terms $\calE(\bff_i), \calE(\bfg_i)$ are quadratically small with tame estimates. This issue also appeared in the work \cite{Damjanovic_Katok10}, but the new difficulty here is that we need to obtain tame splitting in the fiber direction as well as in the base direction. This requires delicate analysis, and our arguments rely on a specific and explicit construction, see Section \ref{Section_tamesplit}. To make up for the fixed loss of regularity,  a standard treatment in the KAM method is that one solves the linearized equations modified by smoothing operators in place of the original linearized equations at each iterative step. To solve problem (III), the intersection property enters into the picture and causes the average terms to be of higher order. See subsection \ref{subsection_induclem}.  Eventually, the iteration process is set and carried out in  subection \ref{subsection_KAMscheme}, which proves Theorem \ref{Element_Thm1}. 

Finally, we go back to prove Theorem \ref{MainThm_0} which is concerned with the perturbation $\tilde\alpha$ of the action $\alpha=\langle \cT_{A_1,\tau_1}, \cT_{A_2,\tau_2}\rangle$. As will be shown in Section \ref{Section_proofMainResult}, the construction of the conjugacy between $\tilde\alpha$ and $\alpha$ mainly consists of two parts: one is produced by the KAM scheme, and the other is obtained by solving a cohomology equation over periodic diffeomorphisms. In general, using only the KAM scheme does not produce the exact conjugacy conjugating $\tilde\alpha$ to $\alpha$. It only produces a conjugacy that conjugates the subgroup $\tilde\alpha\big|_{\mathbf{m}\ZZ+\mathbf{n}\ZZ}$ to  $\alpha\big|_{\mathbf{m}\ZZ+\mathbf{n}\ZZ}$, see Part 1 in the proof of Theorem \ref{MainThm_0}. To solve this issue we need one more step (based on Lemma \ref{lem_periodid} and the commutation relation) to construct a diffeomorphism conjugating the whole action $\tilde\alpha$ to $\alpha$, see Part 2 in the proof of Theorem \ref{MainThm_0}. Theorems \ref{MainThm_1}--\ref{MainThm_2} are obtained as corollaries of Theorem \ref{MainThm_0}.

\subsection{Further discussion} 

Although some of the statements are true in a more general setting, we state them in this paper only for the situation on the torus.

One may wonder if the strategy in this paper could give results for more general actions of similar kind. General question could be: if one can show local rigidity via KAM method for an action on some base manifold,  under which conditions can the KAM method be applied to classify perturbations of fiber bundle extensions of such actions?  More concretely, we may state the following problem.
\begin{Problem}
Let $M$ be a compact nilmanifold and $\rho_0: \ZZ^2\to \textup{Aut}(M)$ be a $\ZZ^2$ action of higher rank.  We consider the extension $\rho=\rho_0\times id_{\TT^s}$ of $\rho_0$ on the bundle $M\times \TT^s$. Suppose that $\rho_0$ is locally rigid via KAM approach. Then, for any smooth action $\tilde\rho:\ZZ^2\to \textup{Diff}^\infty(M\times\TT^s)$ which is sufficiently close to $\rho$ and satisfies the intersection property, is $\tilde\rho$  $C^\infty$-conjugate to $\rho$?
\end{Problem}

That $\rho_0$ is locally rigid via KAM approach means that there is a tame splitting on the base. However, to use the KAM approach for the extended action $\rho=\rho_0\times id_{\TT^s}$, it requires also tame splitting along the fiber direction. In other words, one needs to deal with the following problem. 
\begin{Problem}
Let $C_0^\infty(M\times\TT^s,\RR^s)$ be the space of all smooth functions $u(x,y): M\times\TT^s\to \RR^s$ which satisfy $[u](y):=\int_{M} u(x,y)\, dx=0$. Consider two smooth tame linear operators 
\begin{align*}
d^0: C_0^\infty(M\times\TT^s,\RR^s)&\longrightarrow C_0^\infty(M\times\TT^s,\RR^s)\times C_0^\infty(M\times\TT^s,\RR^s)\\
       u        &\longmapsto (u\circ \rho(\mathbf{e}_1)-u,~u\circ \rho(\mathbf{e}_2)-u)
\end{align*}
\begin{align*}
d^1: C_0^\infty(M\times\TT^s,\RR^s)\times C_0^\infty(M\times\TT^s,\RR^s)&\longrightarrow C_0^\infty(M\times\TT^s,\RR^s)\\
       (u,v)        &\longmapsto u\circ \rho(\mathbf{e}_2)-u-(v\circ \rho(\mathbf{e}_1)-v)
\end{align*}
Does the following exact sequence 
\begin{equation*}
	 C_0^\infty(M\times\TT^s,\RR^s) \xrightarrow{~ d^0~}  C_0^\infty(M\times\TT^s,\RR^s)\times C_0^\infty(M\times\TT^s,\RR^s) \xrightarrow{~d^1~}  C_0^\infty(M\times\TT^s,\RR^s) 
\end{equation*} 
admit a tame splitting? 
\end{Problem}

In the case of $M=\TT^d$, we show in Section \ref{Section_tamesplit} that the tame splitting exists. But we have to say that our proof relies on a specific and explicit construction via Fourier analysis, which allows us to obtain estimates for derivatives along the fiber direction as well the base direction. Nevertheless,  this concrete construction has not yet been generalized to the general compact nilmanifold. This remains a deep and open problem. It may be helpful to use the exponential mixing tool developed in \cite{Gorodnik_Spatzier2015}.

\subsection{Organization of the paper}
This paper is organized as follows. Section \ref{Section_prelim} reviews some basic facts and properties. In Section \ref{Section_actionofproducttype} we consider an action of product type, and state a corresponding local rigidity result (Theorem \ref{Element_Thm1}) for such actions. It plays a crucial role in the proof of Theorem \ref{MainThm_0}. Section \ref{Section_conjugequation}--Section \ref{Section_tamesplit} mainly include obtaining tame solutions for the cohomological equations and constructing tame splitting as well. In Section \ref{Section_KAMscheme}, we prove Theorem \ref{Element_Thm1} by using the KAM scheme.   Theorem \ref{MainThm_0},  Theorem \ref{MainThm_1} and Theorem \ref{MainThm_2} are finally proved in  Section \ref{Section_proofMainResult}.

\subsection{Notation}
In this paper, for a smooth function $f$ we use $\|f\|_{C^r}$ to denote its $C^r$ norm with $r>0$. For smooth functions $f$ and $g$,  we write $\|f\|_{C^r}\leq  C_r \|g\|_{C^r}$ if there exists a sequence of constants $C_r>0$ depending on the regularity $r$ such that these inequalities hold. Accordingly, by $\|f\|_{C^r}\leq  C \|g\|_{C^r}$ we mean that $C$ is a constant which does not vary with $r$. We also write $\|f\|_{C^r}\leq C_{r,s} \|g\|_{C^s}$ in order to stress that the constants depend on both $r$ and $s$.

\section{Preliminaries}\label{Section_prelim}
\subsection{Ergodic toral automorphisms and partial hyperbolicity}

An automorphism of the torus $\TT^d=\RR^d/\ZZ^d$ is determined by a $d\times d$ matrix $A\in \textup{GL}(d,\ZZ)$ with integer entries and determinant $\pm 1$. In this paper, by a slight abuse of notation, we use $A$ to denote both the matrix $A$ and the induced automorphism of $\TT^d$.  The dual to $A$ is the automorphisms $A^*:\ZZ^d\to\ZZ^d$  given by the matrix $A^*=(A^T)^{-1}$. In particular, the Fourier coefficients of any function $\phi\in C^(\TT^d,\RR)$ satisfy: $\widehat{(\phi\circ A)}_n=\hat\phi_{A^*n}$,  $\forall ~ n \in \ZZ^d$.

The following properties are classical, see for instance \cite{Katok_Nitica_2011}. 

\begin{Lem}\label{Lem_charergod}
	(i) An automorphism of $\TT^d$ induced by a matrix $A$ is ergodic if and only if none of the eigenvalues of  $A$ is a root of unity.\\
	(ii) An automorphism of $\TT^d$ induced by a matrix $A$ is ergodic if and only if for any  $n\in\ZZ^d\setminus\{0\}$, the dual orbit $\mathcal{O}(n):=\{(A^*)^i n~:~ i\in\ZZ \}$ is an infinite sequence.\\
	(iii) Any ergodic automorphism  of $\TT^d$ is partially hyperbolic.
\end{Lem}

We infer from Lemma \ref{Lem_charergod} (i) that if $A$ is ergodic, then the automorphism of $\TT^d$ induced by $A^*=(A^T)^{-1}$ is also ergodic. In addition,  by the partial hyperbolicity it has an invariant splitting of the tangent space        
\[\RR^d=E^u(A^*)\oplus E^c(A^*)\oplus E^s(A^*), \] 
and there are a rate $\rho>1$ and a positive constant $C$ such that
\begin{equation}\label{parhyp_split}
\begin{aligned}
    v\in E^u(A^*) ~&\Longleftrightarrow ~\|(A^*)^iv\|\geq C\rho^i\|v\|,  \quad \textup{for all~} i\geq 0, \\
    v\in E^s(A^*) ~& \Longleftrightarrow ~\|(A^*)^iv\|\geq C\rho^{-i}\|v\|,  \quad \textup{for all~} i\leq 0,\\ 
    v\in E^c(A^*) ~& \Longleftrightarrow ~\|(A^*)^iv\|\geq C\frac{\|v\|}{(1+|i|)^d}, \quad \textup{for all~}i\in\ZZ,
\end{aligned}
\end{equation}
 Here,  the superscripts $c, u$ and $s$ stand for ``center'', ``unstable'' and ``stable'', respectively. $A^*$ expands $E^u(A^*)$ (resp. contracts $E^s(A^*)$) with the expansion (resp. contraction) rate being at least $\rho$.
 
For $n\in\ZZ^d$ we can write $n=\pi_u(n)+\pi_s(n)+\pi_c(n),$ where $\pi_u(n), \pi_s(n)$ and $\pi_c(n)$ are the projections of $n$ to the subspaces $E^u(A^*)$, $E^s(A^*)$ and $E^c(A^*)$, respectively. In this paper we say
\begin{itemize}
	\item $n$ is \emph{mostly in} $E^u(A^*)$ and will write $n\hookrightarrow E^u(A^*)$, if  $\|\pi_u(n)\|=\max\limits_{\iota=u,c,s}\|\pi_\iota (n) \|$;
	\item $n$ is \emph{mostly in} $E^s(A^*)$ and will write $n\hookrightarrow E^s(A^*)$, if  $\|\pi_s(n)\|=\max\limits_{\iota=u,c,s}\|\pi_\iota (n) \|$;
	\item $n$ is \emph{mostly in} $E^c(A^*)$ and will write $n\hookrightarrow E^c(A^*)$, if  $\|\pi_c(n)\|=\max\limits_{\iota=u,c,s}\|\pi_\iota (n) \|$.
\end{itemize} 
Obviously, if $n\hookrightarrow E^\iota(A^*)$ with $\iota\in \{u,s,c\}$, then
\begin{equation}\label{ineq_mostlyin}
	\frac{1}{3}\|n\|\leq \|\pi_\iota(n)\|\leq \|n\|.
\end{equation}

The following result comes from the Katznelson lemma \cite{Katznelson_1971}. See also \cite[Remark 5]{Damjanovic_Katok10}.
\begin{Lem}\label{Lem_Katznelson}
Let $M:\TT^d\to\TT^d$ be an ergodic automorphism. Let $V=E^s(M)\oplus E^c(M)$ be the linear subspace in $\RR^d$ spanned by the contracting  and  neutral spaces, then $V\cap \ZZ^d=\{0\}$ and there is a constant $\gamma>0$ such that for any nonzero $n\in \ZZ^d$,
\[\qquad\|\pi_u (n)\|\geq \gamma \,\|n\|^{-d},\] 
where $\pi_u$ is the projection to the expanding space $E^u(M)$, and $\|\cdot\|$ is the Euclidean norm.
\end{Lem}

\subsection{Higher rank actions on tori}
Let us consider higher rank actions by automorphisms of $\TT^d$. 
Recall that a smooth $\ZZ^k$ action $\rho$ by automorphisms of $\TT^d$ is given by a group morphism  $\rho:\mathbf{n}\to  \rho(\mathbf{n})$ from  $\ZZ^k$  into the group $\textup{Aut}(\TT^d)$ of automorphisms of $\TT^d$.  It is a classical result that $\rho$ is higher rank if and only if $\rho(\ZZ^k)$ contains a subgroup isomorphic to $\ZZ^2$ such that all non-trivial elements in this subgroup are ergodic automorphisms, cf. \cite{Starkov1999}.

In particular, in the case of $\ZZ^2$ actions we can say a little more. The following result is elementary, and we give a proof for the reader's convenience.    
\begin{Lem}\label{lem_equihighrank}
	Let $\rho_0=\langle A_1, A_2\rangle=\{A_1^l A_2^k: (l,k)\in \ZZ^2\}$ be a $\ZZ^2$ action generated by automorphisms $A_1$ and $A_2$ on $\TT^d$. If $\rho_0$ is higher rank (i.e., has no rank-one factors), then for any nonzero $(l,k)\in\ZZ^2\setminus\{\mathbf{0}\}$,  $A_1^l A_2^k$ is ergodic on $\TT^d$.
\end{Lem}

\begin{proof}
By assumption, there exists a subgroup $H=\langle \rho_0(\mathbf{i}),\rho_0(\mathbf{j})\rangle$, isomorphic to $\ZZ^2$, such that every non-trivial element in $H$ is ergodic. Here, $\mathbf{i}$ and $\mathbf{j}$ are integer vectors in $\ZZ^2$. 

Now, assume by contradiction that for some $\mathbf{k}\in \ZZ^2\setminus\{0\}$, $\rho_0(\mathbf{k})$ is not ergodic. By Lemma \ref{Lem_charergod}, a toral automorphism is ergodic if and only if none of its eigenvalues is a root of unity. As a consequence,  $\rho_0(n\mathbf{k})$ are non-ergodic for all $n\in\ZZ$.  On the other hand, we observe that 
the subgroups $\langle \mathbf{k} \rangle$ and $\langle \mathbf{i},  \mathbf{j}\rangle$ are, respectively, rank-one  and  rank-two in $\ZZ^2$, so  the intersection between $\langle \mathbf{k} \rangle$ and $\langle \mathbf{i},  \mathbf{j}\rangle$ must be non-trivial and rank-one. Thus, $H$ contains non-identity elements that are non-ergodic. This is a contradiction.
\end{proof}

\begin{Lem}\label{lem_Abexp}
If the $\ZZ^2$ action $\langle A, B\rangle$ generated by automorphisms $A$ and $B$ on $\TT^d$ is higher rank, then there exist constants $\kappa_0>0$ and $C>0$ such that for every non-zero $n\in \ZZ^d$, we have
\begin{equation*}
	\|(A^*)^l(B^*)^k n\|\geq C e^{|(l,k)|\kappa_0}\,\|n\|^{-d},\qquad \textup{for all}~(l,k)\in \ZZ^2.
\end{equation*}
Here, the norm $|(l,k)|:=\max\{|l|, |k|\}$.
\end{Lem}
We refer to \cite{Katok_Katok2005} for the proof.

\subsection{Fr\'echet spaces and tame linear maps}
A Fr\'echet space $X$ is said to be \emph{graded} if the topology is defined by a family  of semi-norms $\{\|\cdot\|_r\}_r$ satisfying $\|x\|_r\leq \|x\|_{r+k}$ for every $x\in X$, and $r, k\geq 0$. For example, the space $C^\infty(\TT^n,\RR)$ with the topology given by the usual $C^r$ norms $\|g\|_r=\max_{0\leq|j|\leq r}\sup_{z\in\TT^n}|\partial^j g(z)|$, $r\in \NN$ is a graded Fr\'echet space. A map $L: U\to V$ between two graded Fr\'echet spaces $U$ and $V$ is said to be \emph{tame} if there exists a constant $\sigma\geq 0$ such that  for any $u\in U$ and $r\geq 0$, 
\[\|L u\|_r\leq C_r\|u\|_{r+\sigma},\]
where the constants $C_r$ may depend on $r$.

Our KAM strategy needs the following classical result (see for instance \cite{Zeh_generalized1,Sal-Zeh} for its proof). It will be used to compensate for the loss of regularity during the KAM iteration.

\begin{Lem}\label{Lem_trun}
There exists a family of linear smoothing operators $\{\rS_N\}_{N\geq 0}$ from $C^\infty(\TT^n,\RR)$ into itself, such that for every $\psi\in  C^\infty(\TT^n,\RR)$, one has $\lim_{N\to\infty}\|\psi-\rS_N \psi\|_{C^0}=0$, and 
\begin{align}\label{trun_est0}
\|\rS_N \psi\|_{C^l}&\leq C_{k,l} N^{l-k}\|\psi\|_{C^k}\qquad \text{for~}  l\geq k,
\end{align}
and for the linear operator $\rR_N\overset{\textup{def}}=id-\rS_N$, it satisfies 
\begin{align} \label{trun_est1}
		\| \rR_N  \psi\|_{C^k}&\leq C_{k,l} \frac{\|\psi\|_{C^l}}{N^{l-k}} \qquad \text{for~}  l\geq k.
\end{align}
Here,  $C_{k,l}>0$ are  constants depending on $k$ and $l$. 
\end{Lem}

In fact, the smoothing operators $\rS_N$ are constructed by convoluting with appropriate kernels decaying  rather fast at infinity. As pointed out in \cite{Zeh_generalized1}, one important consequence of the existence of smoothing operators is the  following interpolation inequalities (Hadamard convexity inequalities).
\begin{Lem}\label{cor_intpest}
	For every  $\psi\in C^\infty(\TT^n,\RR)$ and any $m_1\leq m_2\leq m_3$, $m_2=(1-\lambda) m_1+\lambda m_3$ with $\lambda\in[0,1]$,
	\[\|\psi\|_{C^{m_2}}\leq C_{\lambda,m_1,m_3}\,\|\psi\|_{C^{m_1}}^{1-\lambda}\,\|\psi\|_{C^{m_3}}^{\lambda},\] 
 where the positive constants $C_{\lambda,m_1,m_3}$ depend on $m_1,m_3$ and $\lambda$.
\end{Lem}

We have the following fact on the inverse functions. See for instance \cite{Hamil_1982}.
 \begin{Lem}\label{Apdix_pro1}
 Let $u\in C^\infty(\TT^n, \RR^n)$ and suppose that $\|u\|_{C^1}\leq \frac{1}{4}$. Then, the map induced by $H=id+u:\TT^n\to \TT^n$ is a $C^\infty$ diffeomorphism. Moreover,  the inverse map $H^{-1}$ satisfies 
 \[\|H^{-1}-id\|_{C^r}\leq C_r\,\|u\|_{C^r} \quad\textup{for every~} r\geq 0,\]
where $C_{r}>0$ are constants depending on $r$.
 \end{Lem}

For the composition of two maps, the following estimates hold. 
  \begin{Lem}\label{Apdix_linecont}
 Let  $\psi_1 : B^m\to  B^n $  and  $\psi_2: B^l\to B^m$ be $C^\infty$  maps where $B^\iota\subset \RR^\iota$, $\iota=m, n,l$ are bounded  balls. Suppose that $\|\psi_1\|_{C^1}\leq M$ and $\|\psi_2\|_{C^1}\leq M$ for a constant $M>0$, then the composition $\psi_1\circ \psi_2$ satisfies: for all $r\geq 0$,   
 \begin{align*}
 	\|\psi_1\circ \psi_2\|_{C^r} \leq C_{M,r}\left(1+\|\psi_1\|_{C^r}+\|\psi_2\|_{C^r}\right),
 \end{align*}
where the constants $C_{M,r}$ depend  on $M$ and $r$.
 \end{Lem}
We refer to \cite[Lemma 2.3.4]{Hamil_1982} for its proof.

\section{Partially hyperbolic actions of product type}\label{Section_actionofproducttype}

In this section, we show that our $\ZZ^2$ action $\alpha=\langle \cT_{A_1,\tau_1},\cT_{A_2,\tau_2}\rangle$ can, up to a smooth conjugacy, reduce to an action of product type. The philosophy behind this phenomenon is simple: the higher rank condition on the base space implies that $\tau_i(x)$ is a coboundary with respect to the base map $A_i$, $i=1,2$. More precisely, we can obtain the following result.   
\begin{Pro}\label{Pro_conj_ave}
If $\cT_{A_1,\tau_1}$ commutes with	$\cT_{A_2,\tau_2}$, and $A_1^lA_2^k$ are  ergodic automorphisms on $\TT^d$ for all nonzero $(l,k)\in \ZZ^2$, then there exists a diffeomorphism $\mathfrak{H}\in \textup{Diff}^\infty(\TT^d\times\TT^s)$ which is of the form $\mathfrak{H}(x,y)=(x,y+\phi(x)~\textup{mod}~\ZZ^s)$ with $\phi\in C^\infty(\TT^d,\RR^s)$, such that 
\begin{equation}\label{fbkdufbe}
\mathfrak{H}\circ \cT_{A_1,\tau_1}\circ\mathfrak{H}^{-1}=\cT_{A_1, [\tau_1]},\qquad \mathfrak{H}\circ \cT_{A_2,\tau_2}\circ \mathfrak{H}^{-1}=\cT_{A_2, [\tau_2]},
\end{equation}
where $\mathfrak{H}^{-1}(x,y)=(x,y-\phi(x)~\textup{mod}~\ZZ^s)$, and $[\tau_i]=\int_{\TT^d}\tau_i(x)\, dx$, $i=1,2$.
\end{Pro}

\begin{proof}
Recall that the functions $\tau_1, \tau_2\in C^\infty(\TT^d,\RR^s)$. It is easy to check that $\cT_{A_i,\tau_i}$, $i=1,2$, are conjugate to $\cT_{A_i,[\tau_i]}$ via a common $C^\infty$ conjugacy of the form $\mathfrak{H}(x,y)=(x, y+\phi(x)~\textup{mod}~\ZZ^s)$ if and only if the smooth function $\phi:\TT^d\to \RR^s$ solves the following two cohomological equations
\begin{equation}\label{comsolt}
	\phi(A_1x)-\phi(x)=-\tau_1(x)+[\tau_1],\qquad 	\phi(A_2x)-\phi(x)=-\tau_2(x)+[\tau_2]. 
\end{equation}

Thus, to complete the proof we only need to show that \eqref{comsolt} admits a smooth solution. By the commutation relation $\cT_{A_1,\tau_1}\circ \cT_{A_2,\tau_2}=\cT_{A_2,\tau_2}\circ\cT_{A_1,\tau_1}$, it is direct to see that $A_1$ commutes with $A_2$ and       
\begin{equation}\label{vnnaeb}
	\tau_1(A_2 x)-\tau_1(x)=\tau_2(A_1 x)-\tau_2(x).
\end{equation}
This, combined with the assumption that $A_1^lA_2^k$ are ergodic for all nonzero $(l,k)\in \ZZ^2$, can ensure the existence of $C^\infty$ solutions of equation \eqref{comsolt}. The proof will be included in Lemma \ref{Lem_base_tame2} via Fourier analysis, which provides tame estimates on the solutions as well.	
\end{proof}

\begin{Rem}[Nilmanifold case]
We point out that Proposition \ref{Pro_conj_ave} still holds when the base maps are automorphisms on compact nilmanifolds. The proof requires the use of exponential mixing of the action by automorphisms of nilmanifolds which does not follow easily from Fourier analysis. We refer to the work \cite{Gorodnik_Spatzier2015} by Gorodnik and Spatzier.
\end{Rem}

Observe that $\cT_{A_i,[\tau_i]}$, $i=1,2$, are actually maps of \textbf{product type} since
\[\cT_{A_i,[\tau_i]}=A_i\times R_{[\tau_i]}:\TT^d\times\TT^s\to \TT^d\times\TT^s\] with  $R_{[\tau_i]}$ a translation map on $\TT^s$. In particular, in the case of \textbf{rational} $[\tau_i]$, Proposition \ref{Pro_conj_ave} immediately implies the following result.

\begin{Cor}\label{Cor_rationaltoid}
Let $\mathfrak{H}$ be the conjugacy obtained in Proposition \ref{Pro_conj_ave}. If $[\tau_i]$, $i=1,2$ are both rational, i.e. $[\tau_i]\in \QQ^s$, then for any $q_i\in \ZZ$ satisfying $q_i\, [\tau_i]\in \ZZ^s$, the $q_i$-fold composition $\cT_{A_i, \tau_i}^{q_i}:=\cT_{A_i, \tau_i}\circ \cdots\circ \cT_{A_i, \tau_i}$ is $C^\infty$-conjugate to $A_i^{q_i}\times id_{\TT^s}$, that is 
	\begin{equation*}
	\mathfrak{H}\circ\cT_{A_i, \tau_i}^{q_i}\circ\mathfrak{H}^{-1}=A_i^{q_i}\times id_{\TT^s},
	\end{equation*}
	for each $i=1,2$.
\end{Cor}
\begin{proof} 
	By \eqref{fbkdufbe} one has $\mathfrak{H}\circ \cT_{A_i, \tau_i}^{q_i}\circ \mathfrak{H}^{-1}=\cT_{A_i^{q_i}, q_i [\tau_i]}$, and thus $\cT_{A_i^{q_i}, q_i [\tau_i]}=\cT_{A_i^{q_i}, 0}=A_i^{q_i}\times id_{\TT^s}$. 
\end{proof}

Consequently, Corollary \ref{Cor_rationaltoid} leads us to discover the rigidity phenomenon of a $\ZZ^2$ action generated by two commuting automorphisms of the form $\cT_{A,0}$ and $\cT_{B,0}$.  They are maps of product type: $\cT_{A,0}=A\times id$, $\cT_{B,0}=B\times id:  \TT^d\times\TT^s\longrightarrow \TT^d\times\TT^s$
\begin{equation}\label{form_FFGG}
\begin{aligned}
	\TA(x,y)=(Ax,y),&\qquad \TB(x,y)=(Bx,y).
\end{aligned}
\end{equation}

Let us now use the following equivalent form of the higher rank condition on the action $\langle A, B\rangle$:
\begin{enumerate}
\renewcommand{\labelenumi}{\theenumi}
\renewcommand{\theenumi}{\bf(HR)}
\makeatletter
\makeatother
    \item  \label{condHR}  
    \begin{equation*}
A^l B^k \textup{~is ergodic on~} \TT^d \textup{~for any nonzero~} (l,k)\in\ZZ^2.
    \end{equation*}
\end{enumerate}
Such $A$ and $B$ are called ergodic generators. Then we have the following result.

\begin{The}\label{Element_Thm1}
Let the action $\langle A, B\rangle$ on the base $\TT^d$ satisfies condition \ref{condHR}. Then, there exist  $\vep_0=\vep_0(A, B)>0$ and integer $\mu_0=\mu_0(A, B)$ such that: given any smooth $\ZZ^2$ action $\langle\FF,  \GG\rangle$ on $\TT^d\times\TT^s$, if   $\FF$ and $\GG$ satisfy condition \ref{condIP} and 
\begin{align}\label{qqofnakd}
	\|\FF-\TA\|_{C^{\mu_0}}<\vep_0,\qquad \|\GG-\TB\|_{C^{\mu_0}}<\vep_0,
\end{align}
then $\langle\FF, \GG\rangle$ is $C^\infty$-conjugate to $\langle\TA, \TB\rangle$.
\end{The}  

\begin{Rem}
Let us say a little more on the smallness condition \eqref{qqofnakd}. In fact it suffices to require
\begin{align*}
	\|\FF-\TA\|_{C^0}< \vep,\quad \|\GG-\TB\|_{C^0}<\vep,\qquad 	\|\FF-\TA\|_{C^{\mu_0}}<\vep^{-\frac{3}{4}},\quad \|\GG-\TB\|_{C^{\mu_0}}<\vep^{-\frac{3}{4}}
\end{align*} 
with $\vep$ suitably small.	Through the interpolation estimates, this is enough for the convergence of our KAM scheme. See Lemma \ref{Lem_induc_ineq}.
\end{Rem}

The intersection property \ref{condIP} imposed on $\FF$ and $\GG$ is necessary, otherwise the above result may fail. For instance, consider $\FF=(Ax, y+c)$ and $\GG=(Bx, y+c)$ with $c\neq 0$ a constant vector being arbitrarily small, we find that $\FF$ and $\GG$ cannot be conjugate to $\TA$ and $\TB$. 

On the other hand, the unperturbed maps $\TA$ and $\TB$ indeed satisfy condition \ref{condIP}. In fact, for any $d$-dimensional subtorus $\Gamma\subset \TT^d\times\TT^s$ that is diffeomorphic and $C^1$-close to $\TT^d\times\{y_0\}$, we can write it in the form $\Gamma=\{(x,y)~:~ y=y_0+\psi(x), x\in\TT^d\}$ with $\psi\in C^1(\TT^d,\RR^s)$ suitably small, then the point 
$(0,y_0+\psi(0))$ is exactly a fixed point of $\TA$, which implies $\TA(\Gamma)\cap \Gamma\neq \emptyset$. This is why the intersection property holds for $\TA$. The same is true for $\TB$.

Theorem \ref{Element_Thm1} plays an essential role in proving Theorem \ref{MainThm_0}. In fact, the main task of Sections \ref{Section_conjugequation}--\ref{Section_KAMscheme} is to prove Theorem \ref{Element_Thm1}. The proof is based on the KAM approach. 

\section{The linearized conjugacy equations}\label{Section_conjugequation}

\subsection{Cohomological equations over non-ergodic partially hyperbolic systems}  
In this subsection we will produce the corresponding cohomological equations over a $\ZZ^2$ action by toral automorphisms.  To prove Theorem \ref{Element_Thm1} one needs to find a smooth near-identity diffeomorphism $U$ such that 
\begin{equation}\label{eqnsano}
	U\circ \FF= \TA\circ U,\qquad U\circ \GG=\TB\circ U.
\end{equation}

We introduce in place of $U$ its inverse $H=U^{-1}$ and then write \eqref{eqnsano} in the form
\begin{equation}\label{nciaoqm}
	\FF\circ H= H\circ \TA ,\qquad \GG\circ H=H\circ\TB.
\end{equation}
Writing $H=id+\bfh$ with $\bfh=(\bfh_1,\bfh_2)$, $\bfh_1(x,y)\in C^\infty(\TT^d\times\TT^s,\RR^d)$ and $\bfh_2(x,y)\in C^\infty(\TT^d\times\TT^s,\RR^s)$, and $\FF=\TA+\bff$ and $\GG=\TB+\bfg$,  then \eqref{nciaoqm} reduces to 
\begin{align*}
	\bfh_1\circ \TA-A\bfh_1=\bff_1\circ H,\qquad \bfh_2\circ \TA-\bfh_2=\bff_2\circ H
\end{align*}
and
\begin{align*}
	\bfh_1\circ \TB-B\bfh_1=\bfg_1\circ H,\qquad \bfh_2\circ \TB-\bfh_2=\bfg_2\circ H
\end{align*}
where $\bff_1(x,y), \bfg_1(x,y)\in C^\infty(\TT^d\times\TT^s,\RR^d)$ and $\bff_2(x,y), \bfg_2(x,y)\in C^\infty(\TT^d\times\TT^s,\RR^s)$.
Further, the corresponding linearized equations are 
\begin{equation}\label{twis_LE}
\begin{aligned}
	&\bfh_1(Ax, y)-A\bfh_1(x,y)=\bff_1(x,y),\\
	&\bfh_1(Bx,y)-B\bfh_1(x,y)=\bfg_1(x,y).
\end{aligned}	
\end{equation}
and 
\begin{equation}\label{untwis_LE}
\begin{aligned}
 & \bfh_2 (Ax,y)-\bfh_2(x,y)=\bff_2(x,y),\\
	& \bfh_2(Bx,y)-\bfh_2(x,y)=\bfg_2(x,y).\
\end{aligned}
\end{equation}
Each equation in \eqref{twis_LE} is called a \textit{twisted cohomological equation}, and  each equation in \eqref{untwis_LE} is called an \emph{untwisted cohomological equation}. 

 We point out the following equivalence relation.
\begin{Pro}
Equations \eqref{twis_LE} are solvable in the $C^\infty$ category if and only if the operator 
\begin{equation}\label{fas1}
	\cL_1(\bff_1, \bfg_1)\overset{\textup{def}}=\Big(\bff_1(Bx,y)-B\bff_1(x,y)\Big)-\Big(
	  \bfg_1(Ax,y)-A\bfg_1(x,y)\Big)=0.
\end{equation}
Equations \eqref{untwis_LE} are solvable in the $C^\infty$ category if and only if the operator 
\begin{equation}\label{fas2}
	\cL_2(\bff_2, \bfg_2)\overset{\textup{def}}=\Big(\bff_2(Bx,y)-\bff_2(x,y)\Big)-\Big(\bfg_2(Ax,y)-\bfg_2(x,y)\Big)=0,
\end{equation}
and $\int_{\TT^d}\bff_2(x,y)\,dx=\int_{\TT^d}\bfg_2(x,y)\,dx=0$.
\end{Pro}
\begin{proof}
	It follows directly from Propositions \ref{Pro_LRS0}--\ref{Pro_LRS1111}, which will be shown in Section \ref{Section_Smoothdependparamt}.
\end{proof}

However, we have to say $\cL_1(\bff_1, \bfg_1)\neq 0$ and $\cL_2(\bff_2, \bfg_2)\neq 0$ in general. Instead, they are actually quadratic, see Lemma \ref{Lem_comm_est} below. 

\begin{Lem}\label{Lem_comm_est}
	For the commuting maps $\FF=\TA+\bff$ and $\GG=\TB+\bfg$, we have
	\begin{equation}\label{twoineq_L1L2}
		\| \cL_1(\bff_1,\bfg_1)  \|_{C^r}\leq C_r \|\bff, \bfg\|_{C^{r+1}} \|\bff, \bfg\|_{C^r} ,\qquad \| \cL_2(\bff_2,\bfg_2)  \|_{C^r}\leq C_r \|\bff, \bfg\|_{C^{r+1}} \|\bff, \bfg\|_{C^r}
	\end{equation}
\end{Lem}
\begin{proof}
	By the commutation relation  $\FF\circ \GG=\GG\circ\FF$,  one has
	\begin{align*}
		A\bfg_1+\bff_1\circ\GG=B\bff_1+\bfg_1\circ\FF,\qquad
		\bfg_2+\bff_2\circ\GG=\bff_2+\bfg_2\circ\FF.
	\end{align*}
This implies 
\begin{align}
		\cL_1(\bff_1,\bfg_1) =&\bfg_1\circ\FF-\bfg_1\circ\TA-(\bff_1\circ\GG-\bff_1\circ\TB)
		= \int_0^1 D\bfg_1(\TA+t\bff)\,\bff-D\bff_1(\TB+t\bfg)\,\bfg  dt  
\label{dgdfintgral}
\end{align}
and
\begin{align}
		\cL_2(\bff_2,\bfg_2) =&\bfg_2\circ\FF-\bfg_2\circ\TA-(\bff_2\circ\GG-\bff_2\circ\TB)
		= \int_0^1 D\bfg_2(\TA+t\bff)\,\bff-D\bff_2(\TB+t\bfg)\,\bfg  dt  
\label{dgdfintgral111}
\end{align}
so $\| \cL_1(\bff_1,\bfg_1)  \|_{C^0}\leq C \|\bff, \bfg\|_{C^1} \|\bff, \bfg\|_{C^0}$ and $\| \cL_2(\bff_2,\bfg_2)  \|_{C^0}\leq C \|\bff, \bfg\|_{C^1} \|\bff, \bfg\|_{C^0}$.
This verifies \eqref{twoineq_L1L2}	for $r=0$.  Based on \eqref{dgdfintgral}-\eqref{dgdfintgral111}, the $C^r$-norm estimates follow similarly as in \cite[Lemma 4.7]{Damjanovic_Katok10}.
\end{proof}
 
In view of the quadratic estimates in Lemma \ref{Lem_comm_est}, we can construct approximate solutions of \eqref{twis_LE}--\eqref{untwis_LE} up to errors of higher order. This will be done in Section \ref{Section_tamesplit} and subsection \ref{subsection_induclem}, and it plays an important role in our KAM scheme.

\subsection{Cohomological equations over the base map}\label{subsec_cohomeq_base}
As a warm-up, we first investigate the cohomological equations where all functions involved do not depend on the fiber variables $y$. The results stated here will be used as a ``black box'' for a more general situation discussed in Section \ref{Section_Smoothdependparamt}. In the sequel,  $A$ and $B$ are commuting automorphisms of $\TT^d$ satisfying condition \ref{condHR}, see Section \ref{Section_actionofproducttype}.

\subsubsection{Twisted cohomological equations over ergodic automorphisms of $\TT^d$}
For an ergodic automorphism $A:\TT^d\to\TT^d$ which is partially hyperbolic, we consider the following twisted cohomological equation over $A$, with an unknown function $u:\TT^d\to \RR^d$ and a given function $\Phi:\TT^d\to\RR^d$,
\begin{equation*}
 u(Ax)-Au(x)=\Phi(x),\qquad x\in\TT^d. 
 \end{equation*}
Sometimes, for simplicity we use the symbol $\Delta^A$ to denote $ \Delta^Au(x):=u(Ax)-Au(x)$.

For commuting automorphisms $A$ and $B$ of $\TT^d$, we recall the following result.
\begin{Lem}\cite[Lemma 4.4]{Damjanovic_Katok10}\label{Lem_base_tame1}
	For $\Phi(x)\in C^\infty(\TT^d,\RR^d)$, if there exists a function $\Psi(x)\in C^\infty(\TT^d,\RR^d)$ such that
	$L(\Phi, \Psi):=\Delta^B \Phi-\Delta^A \Psi=0$, then the cohomological equation
	\begin{equation}\label{linequphi}
	\Delta^A u=\Phi
	\end{equation}
	has a unique $C^\infty$ solution $u(x)$, which also solves  the equation $\Delta^B u=\Psi.$ Moreover, it satisfies  
\begin{equation}\label{tame_base1}
	\|u\|_{C^r}\leq C_r \|\Phi\|_{C^{r+\sigma_1}},\qquad \textup{for all~} r\geq 0
\end{equation}
for some  $\sigma_1>0$ depending only on the dimension $d$ and the eigenvalues of $A, B$.
The constants $C_r$ depend on $r$.
\end{Lem}
\begin{Rem}[A remark on inequality \eqref{tame_base1}]
In \cite[Lemma 4.4]{Damjanovic_Katok10} it states that the solution satisfies 
$\|u\|_{C^r}\leq C_r \|\Phi, \Psi\|_{C^{r+\sigma_1}}$. However, according to the proof there one can find that $u$ can be controlled by using only $\Phi$. The existence of $\Psi$ is only used to ensure that the obstruction to solving the linear equation \eqref{linequphi} vanishes. One can also understand it from another perspective: the relation $L(\Phi,\Psi)=0$ implies that $\Psi$ is a solution to the linear equation $\Delta^A \Psi= F$, where $F:=\Delta^B \Phi$. Then, $\Psi$ can be controlled linearly by  $F$, and hence  by  $\Phi$.
\end{Rem}

\subsubsection{Untwisted cohomological equations over ergodic automorphisms of $\TT^d$}
For the ergodic automorphism $A$ of $\TT^d$,  we consider the following untwisted cohomological equation over $A$, with an unknown function $u:\TT^d\to\RR^s$ and a given function $\Phi:\TT^d\to\RR^s$,
\begin{equation*}
u(Ax)-u(x)=\Phi(x),\qquad x\in\TT^d. 
 \end{equation*}

Let $C^\infty_0(\TT^d,\RR^s)$ denote the space of all functions $f\in  C^\infty(\TT^d,\RR^s)$ satisfying $\int_{\TT^d} f(x)\,dx=0$. 
\begin{Lem}\label{Lem_base_tame2}
	For $\Phi(x)\in C_0^\infty(\TT^d,\RR^s)$, if there exists  $\Psi(x)\in C_0^\infty(\TT^d,\RR^s)$ such that
	\begin{equation}\label{condvkf}
	\Phi(Bx)-\Phi(x)=\Psi(Ax)-\Psi(x),
	\end{equation}
then the cohomological equation
	\begin{equation}
	u(Ax)-u(x)=\Phi(x)
	\end{equation}
	has a unique $C^\infty$ solution $u$ in $C^\infty_0(\TT^d,\RR^s)$, and it also solves the equation $u(Bx)-u(x)=\Psi(x).$ Moreover, for any $r\geq 0$   
\begin{equation}\label{tame_base2}
	\|u\|_{C^r}\leq C_r \|\Phi\|_{C^{r+d+2}}
\end{equation}
where the constants $C_r$ depend on $r$.
\end{Lem}

We will use the so-called  \textit{higher-rank trick} developed in \cite{Damjanovic_Katok10} to prove it.

\begin{proof}

Condition  \eqref{condvkf} and the cohomological equations  $u(Ax)- u(x)=\Phi(x)$ and $u(Bx)-u(x)=\Psi(x)$  can split into finitely many one-dimensional problems as follows:
\begin{align}\label{akwoqoo}
\theta(Bx)-\theta(x)= \psi(Ax)- \psi(x)
\end{align}
and the equations
\begin{align}
\omega(Ax)-\omega(x)=\theta(x),\qquad   \omega(Bx)- \omega(x)=\psi(x)
\end{align}
where $\theta, \psi\in C_0^\infty(\TT^d, \RR)$, i.e., $[\theta]=[\psi]=0$. 

Passing to Fourier coefficients, \eqref{akwoqoo} becomes
\begin{align*}
	\hat\theta_{B^*n}-\hat\theta_n=\hat\psi_{A^*n}-\hat\psi_n,\quad \qquad\textup{for ~} n\in\ZZ^d\setminus\{0\}.
\end{align*}
By iterating this formula, for each integer $i\in\ZZ$ we obtain $\hat\theta_{(A^*)^iB^*n}-\hat\theta_{(A^*)^{i}n}$ $=$ $\hat\psi_{(A^*)^{i+1}n}-\hat\psi_{(A^*)^{i}n}$. Taking the formal sum over  $i$ we obtain 
\begin{equation}\label{winncnda}
\sum_{i\in\ZZ}\hat\theta_{(A^*)^iB^*n}-	\sum_{i\in\ZZ}\hat\theta_{(A^*)^{i}n}=\sum_{i\in\ZZ}\hat\psi_{(A^*)^{i+1}n}-\sum_{i\in\ZZ}\hat\psi_{(A^*)^{i}n}
\end{equation}
Here, $n$ has nontrivial projections to unstable subspace $E^u(A^*)$ and stable subspace $E^s(A^*)$ since  $A^*$ is partially hyperbolic,  so one can find that all the sums involved in 
\eqref{winncnda} are absolutely convergent when $n\neq 0$ (see  \cite[Lemma 4.3]{Damjanovic_Katok10}). Note that the right-hand side of \eqref{winncnda} equals zero, hence, the left-hand side implies that
\begin{align*}
	\sum_{i\in\ZZ}\hat\theta_{(A^*)^iB^*n}=\sum_{i\in\ZZ}\hat\theta_{(A^*)^{i}n},\qquad\textup{for every~} n\in \ZZ^d\setminus\{0\}.
\end{align*}
By iterating this equation, for each $j\in\ZZ$,
\begin{align}\label{rurdllw}
	\sum_{i\in\ZZ}\hat\theta_{(A^*)^i(B^*)^jn}=\sum_{i\in\ZZ}\hat\theta_{(A^*)^i(B^*)^{j-1}n}=\cdots\cdots=\sum_{i\in\ZZ}\hat\theta_{(A^*)^iB^*n}=\sum_{i\in\ZZ}\hat\theta_{(A^*)^{i}n}.
\end{align}
By condition \ref{condHR} in Section \ref{Section_actionofproducttype} and the ergodicity there,  one can show that	$\sum_{i\in\ZZ}\hat\theta_{(A^*)^i(B^*)^jn}$ converges to zero, as $j\to \infty$. Therefore, \eqref{rurdllw} implies that
\begin{equation}\label{obvani}
	\sum_{i\in\ZZ}\hat \theta_{(A^*)^{i}n}=0, \qquad\textup{for each~} n\in \ZZ^d\setminus\{0\}. 
\end{equation}  

Now, we consider the cohomological equation $\omega(Ax)-\omega(x)=\theta(x)$. Passing to Fourier coefficients, it is equivalent to solving the following equations 
\begin{equation}\label{foucon}
	\hat\omega_{A^*n}-\hat\omega_{n}=\hat\theta_n, \qquad \forall~n\in\ZZ^d.
\end{equation} 
When $n=0$,  we can set $\hat\omega_0=0$ since we have assumed $\hat\theta_0=0$. For $n\in \ZZ^d\setminus\{0\}$, from \eqref{obvani} we see that the obstruction to solving equation  \eqref{foucon} vanishes, so we have $\hat\omega_n= \hat\omega_n^+=\hat\omega_n^-$, where  
\begin{equation}\label{adjsk}
	\hat\omega_n^+=-\sum_{i=0}^\infty\theta_{(A^*)^i n},\qquad  \hat\omega_n^-=\sum_{i= -1}^{-\infty}\theta_{(A^*)^i n}. 
\end{equation}
Consequently, we obtain a formal solution \[\omega=\sum_{n\in\ZZ^d} \hat\omega_n\cdot e^{i2\pi\langle n,x\rangle}\]
for the  equation $\omega(Ax)-\omega(x)=\theta(x)$.  In order to  show that $\omega\in C^\infty$, we need to estimate   $\hat\omega_n$ for every $n\neq 0$.

If $n$ is mostly in $E^u(A^*)$, i.e. $n\hookrightarrow E^u(A^*)$, we  use the form $\hat \omega_n=\hat \omega_n^+$, by \eqref{parhyp_split} and \eqref{ineq_mostlyin},  for any $k\in \ZZ^+$ we have
\begin{align}\label{aanoi1}
	\begin{aligned}
	|\hat \omega_n|	
		\leq \sum_{i\geq 0}\frac{\|\theta\|_{C^k}}{\|(A^*)^in\|^{k}}
		\leq \sum_{i\geq 0}\frac{\|\theta\|_{C^k}}{\|(A^*)^i\,\pi_u(n)\|^{k}}
		\leq \sum_{i\geq 0}\frac{\|\theta\|_{C^k}}{\rho^{ik}\|\pi_u(n)\|^{k}}
		\leq &  \frac{M_k\cdot\|\theta\|_{C^k}}{\|\pi_u(n)\|^{k}}\leq C_k \frac{\|\theta\|_{C^k}}{\|n\|^{k}}
	\end{aligned}
\end{align}
where the expanding rate $\rho>1$ and $C_k=3^k\cdot M_k$. Similarly, if $n$ is mostly in $E^s(A^*)$, i.e. $n\hookrightarrow E^s(A^*)$, we use the form $\hat \omega_n=\hat \omega_n^-$ to obtain 
\begin{equation}\label{aanoi2}
		|\hat \omega_n|\leq C_k \|\theta\|_{C^k}\,\|n\|^{-k},\qquad \forall~k\in\ZZ^+.
\end{equation}

If  $n\hookrightarrow E^c(A^*)$, we use the form $\hat \omega_n=\hat \omega_n^+$. 
By the Katznelson lemma (see Lemma \ref{Lem_Katznelson}) one has $\|\pi_u(n)\|\geq \gamma \|n\|^{-d}$. Then,
\begin{align*}
	\|(A^*)^i n\|\geq\|(A^*)^i\,\pi_u(n)\|\geq C\rho^{i}\|\pi_u(n)\|\geq C \gamma  \rho^{i} \|n\|^{-d}\geq C\gamma  \rho^{i-i_0} \|n\|.
\end{align*}
 for all  $i\geq i_0$ where  $i_0=\left[\frac{(d+1)\ln\|n\|}{\ln \rho}\right]+1$. 
For $0\leq i \leq i_0-1$, by \eqref{parhyp_split} and \eqref{ineq_mostlyin} we have 
\[\|(A^*)^i n\|\geq \|(A^*)^i \pi_c(n)\|\geq C(1+i)^{-d}\|\pi_c(n)\|\geq\frac{C}{3}(1+i)^{-d}\|n\|.  \]
Then it follows that for $k\in\ZZ^+$,
\begin{align}
	|\hat\omega_n| \leq\sum_{i\geq 0}\frac{\|\theta\|_{C^k}}{\|(A^*)^in\|^{k}}
		\leq   C'\frac{\|\theta\|_{C^k}}{\|n\|^{k}}\left(\sum_{i=0}^{i_0-1}(1+i)^{dk}+ \sum_{i=i_0}^{\infty}\rho^{-k(i-i_0)}\right)
		\leq  & C''\|\theta\|_{C^k} \|n\|^{-k}\cdot (i_0^{dk+1}+ c)\nonumber\\
		\leq & C''' \|\theta\|_{C^k} \|n\|^{-k}  \cdot(\ln\|n\|)^{dk+1}\nonumber\\
		\leq &   C_k \|\theta\|_{C^k}\, \|n\|^{-k+1}\label{haomgn2}
\end{align}

Finally, for each $r\geq 0$, using \eqref{aanoi1}--\eqref{haomgn2} we obtain that
\begin{align*}
	\|\omega\|_{C^r}\leq (2\pi)^r \sum_{n\in\ZZ^d\setminus\{0\}} \|n\|^r\cdot|\hat \omega_n|\leq  (2\pi)^r \,C_k \sum_{n\in\ZZ^d\setminus\{0\}}  \frac{\|\theta\|_{C^k}}{\|n\|^{k-r-1}}
\end{align*}
In particular, taking $k=r+d+2$ we obtain \[\|\omega\|_{C^r}\leq C'_r \|\theta\|_{C^{r+d+2}}.\] 
Since this is true for any $r\geq 0$, $\omega\in C_0^\infty(\TT^d,\RR)$, and it solves the equation \[\omega(Ax)-\omega(x)=\theta(x).\] The uniqueness of solutions  follows from the ergodicity of $A$.

Since $A$ and $B$ commute, by \eqref{akwoqoo} we can prove that $\omega$ also solves $\omega(Bx)-\omega(x)=\psi(x)$. The argument is elementary and similar to \cite[Lemma 4.4]{Damjanovic_Katok10}, so we will not repeat here.
\end{proof}

\section{Smooth dependence on multi-dimensional parameters}\label{Section_Smoothdependparamt}

The main goal of this section is to prove Propositions \ref{Pro_LRS0}--\ref{Pro_LRS1111}. In contrast with subection \ref{subsec_cohomeq_base}, we will deal with cohomological equations over a non-ergodic partially hyperbolic automorphism  $\TA=A\times id_{\TT^s}$ with $A$ ergodic on the base $\TT^d$. Based on the results in subection \ref{subsec_cohomeq_base}, we then employ the idea of smooth dependence on parameters to study the solutions of the equations of the form \eqref{twis_LE}--\eqref{untwis_LE} and give their derivative estimates in the fiber direction as well as in the base direction. A similar idea was once used by \cite{delaLave_Marco_Moriyon1986} to study the parameter dependence for the solutions of untwisted cohomological equations over Anosov diffeomorphisms.

Throughout this section, $A:\TT^d\to\TT^d$ and $B:\TT^d\to\TT^d$ are commuting automorphisms and satisfy condition \ref{condHR} given in Section \ref{Section_actionofproducttype}.

\subsection{The twisted case}\label{subsect_thetwistcase}
We first consider the twisted cohomological equations.
Let us introduce the following set in $C^\infty(\TT^d,\RR^d)$.
\begin{equation*}
\mathbb{V}:=\{\phi\in C^\infty(\TT^d,\RR^d)~:\quad \exists~ \psi\in C^\infty(\TT^d,\RR^d) \textup{~such that~} \Delta^B \phi=\Delta^A \psi \} 
\end{equation*}
where the symbols $\Delta^A\,\psi(x)=\psi(Ax)-A \psi(x)$ and  $\Delta^B\,\phi(x)=\phi(Bx)-B \phi(x)$. According to subsection \ref{subsec_cohomeq_base} it is easy to find that $\mathbb{V} $ is exactly the set of all functions $\phi\in C^\infty(\TT^d,\RR^d)$ for which the equation $\Delta^A u=\phi$ admits a smooth solution.

For our purpose, we need to study $s$--dimensional parameters. The set of parameters  will be an open ball  $\cD\subset \RR^s$, and we use $y\in \cD$ to denote the parameter variables.

\begin{Lem}\label{Lem_Mcont}
$\mathbb{V}$ has the following property: 
\begin{enumerate}[(i)]
	\item $\mathbb{V}$ is a linear subspace of $C^\infty(\TT^d, \RR^d)$. 
There is a tame linear operator $\cM: \mathbb{V}\longrightarrow C^\infty(\TT^d,\RR^d)$ which satisfies: for each $\phi\in \mathbb{V}$,
\begin{equation}\label{teqkaslf1}
	\Delta^A \big(\cM(\phi)\big)=\phi \quad \textup{~and~} \quad \|\cM(\phi)\|_{C^r(\TT^d)}\leq C_r \|\phi\|_{C^{r+\sigma_1}(\TT^d)},
\end{equation}
where $\sigma_1>0$ is the same constant given in Lemma \ref{Lem_base_tame1}.
\item For $\xi(x,y)\in C^\infty(\TT^d\times\cD, \RR^d)$, we denote $\xi^y(x):=\xi(x,y)$. 
If $\xi^y\in \mathbb{V}$ for every parameter $y\in \cD$, then the map $y\longmapsto \cM({\xi}^y)$ is continuous, i.e., for any $r\in\NN$,
\begin{align}\label{lemcontin}
	\lim_{y\to a}\| \cM({\xi}^y)- \cM({\xi}^{a})\|_{C^{r}(\TT^d)}=0,\qquad \forall~a\in\cD.
\end{align}
\end{enumerate}
\end{Lem}
\begin{proof}  
(i) The fact that $\mathbb{V}$ is a linear subspace follows readily from the definition. Moreover, \eqref{teqkaslf1}	comes from Lemma \ref{Lem_base_tame1}. Namely, for each $\phi\in \mathbb{V}$, $\cM(\phi)$ is the unique solution of $\Delta^A u=\phi$. 

(ii) Given $r\in \NN$.
Assume by contradiction that there exists some point $a\in\cD$ such that \eqref{lemcontin} fails. Then there would exist a sequence $z_k\to a$ and a number $\delta>0$ such that   
\begin{equation}\label{xnbnfaa}
	\| \cM({\xi}^{z_k})- \cM({\xi}^a)\|_{C^{r}(\TT^d)}> \delta,\quad \textup{for all}~ k.
\end{equation}
On the other hand, by item (i) we see that
\begin{align*}
	\|\cM({\xi}^{z_k})- \cM({\xi}^a)\|_{C^{r+1}(\TT^d)}\leq C_{r+1} \| \xi^{z_k}-\xi^{a}\|_{C^{r+1+\sigma_1}(\TT^d)}.
\end{align*}
Then all the functions $\cM({\xi}^{z_k})$ are uniformly bounded in  the $C^{r+1}$ topology because  $\xi(x,y)\in C^\infty$. Using the Arzel\`a-Ascoli theorem it is not difficult to show that, by taking a subsequence if necessary, $ \cM({\xi}^{z_k})$ converges to some $v$ in  the  $C^r(\TT^d,\RR^d)$ topology. 

However, by continuity we have $\Delta^A v=\xi^a$. Hence the uniqueness of solutions implies 
$v=\cM({\xi}^a)$, which  contradicts \eqref{xnbnfaa}. This finishes the proof.
\end{proof}

For two  smooth functions $f, g: \TT^d\times\cD\to \RR^d $, we define 
\[\cL_1(f, g):=\Big(f(Bx,y)-B f(x,y)\Big)-\Big(g(Ax,y)-A g(x,y)\Big).\]
In the sequel, for a smooth function $\xi(x,y)$ we denote $(\partial_y^\beta\xi)^y(x):=\partial_y^\beta \xi(x,y)$.

\begin{Lem}\label{Lem_xieta}
Suppose that  $\xi, \eta \in C^\infty(\TT^d\times\cD,\RR^d)$ and   $\cL_1(\xi,\eta)=0$. Then, 
\begin{enumerate}[(i)]
\item for any differential operator $\partial^\beta_y$ with the multi-index $\beta\in\NN^s$, we  have 
	$(\partial_y^\beta\xi)^y\in\mathbb{V}$. 
	\item for any  parameter $y\in\cD$ and any index $j=1,\cdots, s$
\begin{align}\label{lemcnt2}
	\lim_{\vep\to 0}\left\|\frac{\cM\left(\xi^{y+\vep \mathbf{e}_j}\right)-\cM(\xi^{y})}{\vep}-\cM \big((\partial_{y_j} \xi)^y\big)\right\|_{C^r(\TT^d)}=0,\quad\textup{for any~} r\in\NN.
\end{align}
where 
$(\mathbf{e}_1, \cdots,\mathbf{e}_j,\cdots,\mathbf{e}_s)$ is the standard orthogonal basis for $\RR^s$.  This also implies that the map $y\longmapsto \cM({\xi}^y)\in C^\infty(\TT^d,\RR^d)$ is of class $C^1$.
\end{enumerate}	
\end{Lem}
\begin{proof}
	(i) 
It follows easily from the  fact  $\cL_1(\partial^\beta_{y} \xi, \partial_{y}^\beta \eta)=\partial_{y}^\alpha\cL_1(\xi,\eta)=0$.

(ii)  We only need to check the case of $j=1$, other cases are similar. If \eqref{lemcnt2} fails  for some $y\in\cD$ and some $r\in\NN$, there would exist a nonzero
sequence $\vep_k\to 0$ in $\RR$ such that
\begin{align}\label{xiyedel}
	\left\|\frac{\cM\left(\xi^{y+\vep_k \mathbf{e}_1}\right)-\cM(\xi^{y})}{\vep_k}-\cM \big((\partial_{y_1} \xi)^y\big)\right\|_{C^r(\TT^d)}> \delta
\end{align}
for some  number $\delta>0$.  By the tame estimate \eqref{teqkaslf1} we have 
\begin{align*}
	\left\|\frac{\cM\left(\xi^{y+\vep_k \mathbf{e}_1}\right)-\cM(\xi^{y})}{\vep_k}-\cM \big((\partial_{y_1} \xi)^y\big)\right\|_{C^{r+1}(\TT^d)}=&\left\|\cM\left(\frac{\xi^{y+\vep_k \mathbf{e}_1}-\xi^{y}}{\vep_k}- (\partial_{y_1} \xi)^y\right)\right\|_{C^{r+1}(\TT^d)}\\
	\leq & C_{r+1+\sigma_1} \left\|\frac{\xi^{y+\vep_k \mathbf{e}_1}-\xi^{y}}{\vep_k}- (\partial_{y_1} \xi)^y\right\|_{C^{r+1+\sigma_1}(\TT^d)}
\end{align*}
Since $\xi\in C^\infty$, the quantity in the last line is uniformly  bounded as $\vep_k\to 0$. This implies 
$\frac{\cM\left(\xi^{y+\vep_k \mathbf{e}_1}\right)-\cM(\xi^{y})}{\vep_k}$ is uniformly bounded in the $C^{r+1}(\TT^d,\RR^d)$ topology. Then, using the Arzel\`a-Ascoli theorem we are able to prove that, up to a subsequence,  $\frac{\cM\left(\xi^{y+\vep_k \mathbf{e}_1}\right)-\cM(\xi^{y})}{\vep_k}$ converges to some $w$ in the  $C^r(\TT^d,\RR^d)$ topology. Due to \eqref{xiyedel}, we get 
 \begin{align}\label{wxidel}
	\left\|w-\cM \big((\partial_{y_1} \xi)^y\big)\right\|_{C^r(\TT^d)}> \delta.
\end{align}

On the other hand, it is direct to see that
\begin{align*}
	\Delta^A \left(\frac{\cM\left(\xi^{y+\vep_k \mathbf{e}_1}\right)-\cM(\xi^{y})}{\vep_k}\right)=\frac{\xi^{y+\vep_k \mathbf{e}_1}-\xi^{y}}{\vep_k}.
\end{align*}
Sending $\vep_k\to 0$ yields $
	\Delta^A w=(\partial_{y_1}\xi)^{y}$. Then the uniqueness of solutions implies that
	$w=\cM \big((\partial_{y_1} \xi)^y\big)$. This contradicts \eqref{wxidel}.
	\end{proof}

Since a smooth function $R(x,y)$ on $\TT^d\times\TT^s$ is also a  smooth function  on $\TT^d\times\RR^s$ that is $\ZZ^s$-periodic in $y$, we can apply the above lemmas to prove the following result.
\begin{Pro}\label{Pro_LRS0}
Suppose that $\cL_1(R, S)=0$ where $R(x,y), S(x,y): \TT^d\times\TT^s\to \RR^d$ are $C^\infty$ functions. Then, there is a unique function $\Omega\in C^\infty(\TT^d\times\TT^s,\RR^d)$ that solves the equations 
		\begin{align}\label{exten_twoeqs}
			\Delta^A \Omega(x,y)= R(x,y),\qquad \Delta^B \Omega(x,y)= S(x,y).
		\end{align}
Here, $\Delta^A \Omega(x,y)=\Omega(Ax,y)-A \Omega(x,y)$ and $\Delta^B \Omega(x,y)=\Omega(Bx,y)-B \Omega(x,y)$. 
Moreover, 
	\begin{align}\label{ext_omegest}
		\| \Omega \|_{C^r}\leq C_r \|R\|_{C^{r+\sigma_1}},
	\end{align}
	where $\sigma_1>0$ is the same constant given in Lemma \ref{Lem_base_tame1}.
\end{Pro}

\begin{proof}

\noindent\textbf{(I) Existence of continuous solutions to \eqref{exten_twoeqs}.}
We use the notation \[R^y(x):=R(x,y),\qquad S^y(x):=S(x,y).\]
The condition $\cL_1(R, S)=0$ implies that $R^y\in \mathbb{V}$, for any $y$.
Then by Lemma \ref{Lem_base_tame1} we obtain a family of solutions $u^y:=\cM( R^y)\in C^\infty(\TT^d,\RR^d)$, with parameter $y$, that solves
 \[ \Delta^A u^y=R^y,\qquad \Delta^B u^y= S^y, \qquad \textup{for each}~ y.  \]
 Thus, we define
  \[\Omega(x,y):=u^y(x)=\cM(R^y)(x).\] 
 $\Omega(x,y)$ is continuous on $\TT^d\times\TT^s$: indeed, for any point $(a,b)$, \begin{align*}
|\Omega(x,y)-\Omega(a,b)|\leq &	|\Omega(x,y)-\Omega(x,b)|+|\Omega(x,b)-\Omega(a,b)|\\
=& |u^y(x)-u^{b}(x)|+|u^{b}(x)-u^{b}(a)|\\
\leq & \|\cM(R^y)-\cM(R^{b})\|_{C^0(\TT^d)}+|\cM(R^{b})(x)-\cM(R^{b})(a)|
\end{align*}
By Lemma \ref{Lem_Mcont} (ii), the last line tends to zero as $(x,y)\longrightarrow (a,b)$. Consequently,  
 $\Omega(x,y)$ is a continuous solution to \eqref{exten_twoeqs}. More precisely, it is smooth in $x$ and continuous in $y$.

\noindent\textbf{(II) $C^1$-regularity. } 
Applying Lemma \ref{Lem_xieta} (ii) to $R^y$ we obtain that  
 for each $j=1,\cdots,s$, the partial derivative $\partial_{y_j}\Omega(x,y)$ exists and
 \[\partial_{y_j}\Omega(x,y)=\cM\big((\partial_{y_j} R)^{y}\big)(x).\]

Now, let us show the continuity of $(x,y)\longmapsto\partial_{y_j}\Omega(x,y)$. Indeed, for any point $(a,b)$, 
\begin{align*}
\left|\partial_{y_j}\Omega(x,y)-\partial_{y_j}\Omega(a,b)\right|\leq &	\left|\partial_{y_j}\Omega(x,y)-\partial_{y_j}\Omega(x, b)\right|+\left|\partial_{y_j}\Omega(x, b)-\partial_{y_j}\Omega(a, b)\right|\\
\leq & \left\|\cM\big((\partial_{y_j} R)^{y}\big)-\cM\big((\partial_{y_j} R)^{b}\big)\right\|_{C^0(\TT^d)}\\
& \qquad\qquad\quad +\left|
\cM\big((\partial_{y_j} R)^{b}\big)(x)-\cM\big((\partial_{y_j} R)^{b}\big)(a)\right|.
\end{align*}
Here, it is evident that $\left|
\cM\big((\partial_{y_j} R)^{b}\big)(x)-\cM\big((\partial_{y_j} R)^{b}\big)(a)\right|$ converges to zero as $x\to a$. Meanwhile, as $(\partial_{y_j} R)^y\in\mathbb{V}$ for all $y$,
using Lemma \ref{Lem_Mcont} (ii) with $\xi=\partial_{y_j} R$ we deduce that
  \[\left\|\cM\big((\partial_{y_j} R)^{y}\big)-\cM\big((\partial_{y_j} R)^{b}\big)\right\|_{C^0(\TT^d)}\longrightarrow 0 \]  as $y\to b$. Thus,  $\partial_{y_j}\Omega(x,y)$ converges to $\partial_{y_j}\Omega(a,b)$ as  $(x,y)\to (a,b)$.

On the other hand, 
to prove  the continuity of  $(x,y)\longmapsto\partial_{x_i}\Omega(x,y)$ with $x_i$ the $i$-th coordinate of $x$, we observe that
\begin{align*}
\left|\partial_{x_i}\Omega(x,y)-\partial_{x_i}\Omega(a,b)\right|\leq &	\left|\partial_{x_i}\Omega(x,y)-\partial_{x_i}\Omega(x, b)\right|+\left|\partial_{x_i}\Omega(x, b)-\partial_{x_i}\Omega(a, b)\right|\\
\leq & \left\|\cM\big(R^{y}\big)-\cM\big(R^{b}\big)\right\|_{C^1(\TT^d)}
 +\left|
\partial_{x_i}\cM\big(R^{b}\big)(x)-\partial_{x_i}\cM\big(R^{b}\big)(a)\right|,
\end{align*}
where the last line converges to zero as $(x,y)\longrightarrow(a,b)$, thanks to Lemma \ref{Lem_Mcont} (ii) and $\cM\big(R^{b}\big)\in C^\infty(\TT^d,\RR^d)$.
Therefore, we conclude that $\Omega(x,y)$ is of class $C^1$.

\noindent\textbf{(III) $C^k$-regularity. } 
The higher regularity can be proved by induction. 
The case of $r=1$ has been proved above.
Suppose that $\Omega(x,y)$ is $C^r$ and 
\[\partial_x^\alpha\partial_y^{\beta}\Omega(x,y)=\partial_x^\alpha \cM\big((\partial^\beta_y R)^y\big)\] 
for any multi-indices $\alpha, \beta$ satisfying $|\alpha|+|\beta|=r$, we will show that  $\Omega\in C^{r+1}$, namely, every partial derivatives of order  $\leq r+1$ exists and is continuous.
Since $\Omega$ is assumed to be $C^r$, one only needs to check that  the  partial derivatives $\partial_{x_i}\partial_x^\alpha\partial_y^{\beta}\Omega(x,y)$ and $\partial_{y_j}\partial_x^\alpha\partial_y^{\beta}\Omega(x,y)$, $|\alpha|+|\beta|=r$, exist and continuous on $\TT^d\times\TT^s$.

 We first claim that the partial derivative $\partial_{y_j}\partial_x^\alpha\partial_y^\beta\Omega(x,y)$, where $|\alpha|+|\beta|=r$, exists and
\begin{equation}\label{bvnsb}
	\partial_{y_j}\partial_x^\alpha\partial_y^\beta\Omega(x,y)=\partial^\alpha_x\cM\big((\partial_{y_j}\partial_y^\beta R)^{y}\big).
\end{equation}
Let $(\mathbf{e}_1, \cdots,\mathbf{e}_j,\cdots,\mathbf{e}_s)$ denote the standard orthogonal basis for $\RR^s$.
 Then, for $\vep\neq 0$,
 \begin{equation}\label{limtepej}
 \begin{aligned}
 &\left|\frac{\partial_x^\alpha\partial_y^\beta\Omega(x,y+\vep\mathbf{e}_j)-\partial_x^\alpha\partial_y^\beta\Omega(x,y)}{\vep}-  \partial^\alpha_x\cM\big((\partial_{y_j}\partial_y^\beta R)^{y}\big)\right|\\  
 &\qquad\qquad\qquad  = \left|\partial_x^\alpha\left(\frac{\partial_y^\beta\Omega(x,y+\vep\mathbf{e}_j)-\partial_y^\beta\Omega(x,y)}{\vep}-\cM\big((\partial_{y_j}\partial_y^\beta R)^{y}\big)\right)\right|\\
&\qquad\qquad\qquad \leq \left\|\frac{\cM\big((\partial_y^\beta R)^{y+\vep\mathbf{e}_j}\big)-\cM\big((\partial_y^\beta R)^{y}\big)}{\vep}-\cM\big((\partial_{y_j}\partial_y^\beta R)^{y}\big)\right\|_{C^{|\alpha|}(\TT^d)}
\end{aligned} 
\end{equation}
Note that $\partial^\beta_y R(x, y)\in C^\infty$ and $\cL_1(\partial^\beta_y R, \partial^\beta_y S)=0$. Then, applying Lemma \ref{Lem_xieta} (ii) with $\xi=\partial^\beta_y R$ we see that the last line of \eqref{limtepej} converges to zero as $\vep\to 0$. Therefore, $\partial_{y_j}\partial_x^\alpha\partial_y^\beta\Omega(x,y)$ exists and \eqref{bvnsb} holds. 

Next, we will show that $(x,y)\mapsto$ $\partial_{y_j}\partial_x^\alpha\partial_y^\beta\Omega(x,y)$ is continuous. 
Indeed,  for any point $(a, b)$, using \eqref{bvnsb} it follows that
\begin{align*}
	&\left|\partial_{y_j}\partial_x^\alpha\partial_y^\beta\Omega(x,y)-\partial_{y_j}\partial_x^\alpha\partial_y^\beta\Omega(a,b)\right|\\
	\leq & \left|\partial_{y_j}\partial_x^\alpha\partial_y^\beta\Omega(x,y)-\partial_{y_j}\partial_x^\alpha\partial_y^\beta\Omega(x,b)\right|+\left|\partial_{y_j}\partial_x^\alpha\partial_y^\beta\Omega(x,b)-\partial_{y_j}\partial_x^\alpha\partial_y^\beta\Omega(a,b)\right|\\
	= &\left|\partial^\alpha_x\cM\big((\partial_{y_j}\partial_y^\beta R)^{y}\big)-\partial^\alpha_x\cM\big((\partial_{y_j}\partial_y^\beta R)^{b}\big)\right|+\left|\partial^\alpha_x\cM\big((\partial_{y_j}\partial_y^\beta R)^{b}\big)(x)-\partial^\alpha_x\cM\big((\partial_{y_j}\partial_y^\beta R)^{b}\big)(a)\right|\\
	\leq & \left\|\cM\big((\partial_{y_j}\partial_y^\beta R)^{y}\big)-\cM\big((\partial_{y_j}\partial_y^\beta R)^{b}\big)\right\|_{C^{|\alpha|}}+\left|\partial^\alpha_x\cM\big((\partial_{y_j}\partial_y^\beta R)^{b}\big)(x)-\partial^\alpha_x\cM\big((\partial_{y_j}\partial_y^\beta R)^{b}\big)(a)\right|.		
\end{align*}
Evidently, the second quantity of the last line converges to zero as $x\to a$. By applying Lemma \ref{Lem_Mcont} (ii) with $\xi=\partial_{y_j}\partial_y^\beta R$, the first quantity of the last line also converges to zero as $y\to b$.
This thus verifies the continuity of $(x,y)\mapsto$ $\partial_{y_j}\partial_x^\alpha\partial_y^\beta\Omega(x,y)$. 

By using similar argument, we can also prove that the partial derivative $\partial_{x_i}\partial_x^\alpha\partial_y^\beta\Omega(x,y)$, $|\alpha|+|\beta|=r$, exists and equals $
	\partial_{x_i}\partial^\alpha_x\cM\big((\partial_y^\beta R)^{y}\big)
$, and the function $(x,y)\mapsto$ $\partial_{x_i}\partial_x^\alpha\partial_y^\beta\Omega(x,y)$ is continuous.  Therefore we can conclude  that $\Omega$ is $C^{r+1}$. 

By what we have shown above, equations \eqref{exten_twoeqs} have a unique $C^\infty$ solution $\Omega$. In addition, for any multi-indices $\alpha, \beta$ with $|\alpha|+|\beta|=r$,  as $\partial_x^\alpha\partial_y^{\beta}\Omega(x,y)=\partial_x^\alpha \cM\big((\partial^\beta_y R)^y\big)$, we can apply Lemma \ref{Lem_Mcont} to obtain that
\begin{align*}
	|\partial_x^\alpha\partial_y^{\beta}\Omega|=|\partial_x^\alpha \cM\big((\partial^\beta_y R)^y\big)|\leq  \|\cM\big((\partial^\beta_y R)^y\big)\|_{C^{|\alpha|}(\TT^d)}\leq C_{|\alpha|} \|(\partial^\beta_y R)^y\|_{C^{|\alpha|+\sigma_1}(\TT^d)} 
	\leq  C_r \|R\|_{C^{r+\sigma_1}},
\end{align*}
which finally verifies estimate \eqref{ext_omegest}.
\end{proof}

\subsection{The untwisted case}
The untwisted cohomological equations can be studied in the same spirit. 
Let 
\[\cL_2(R, S):=\Big(R(Bx,y)-R(x,y)\Big)-\Big(S(Ax,y)-S(x,y)\Big),\] we have the following property. 
\begin{Pro}\label{Pro_LRS1111}
Suppose that $\cL_2(R, S)=0$, where $R(x,y), S(x,y): \TT^d\times\TT^s\to \RR^s$ are $C^\infty$ functions and the averages over the base $\TT^d$ vanish: $\int_{\TT^d} R(x,y) \,dx=\int_{\TT^d} S(x,y)\, dx=0$. Then, 		\begin{align}
			\Omega(Ax,y)-\Omega(x,y)= R(x,y),\qquad  \Omega(Bx,y)-\Omega(x,y)= S(x,y)
		\end{align}
	have a unique solution $\Omega\in C^\infty(\TT^d\times\TT^s,\RR^s) $ satisfying $\int_{\TT^d} \Omega(x,y)\,dx=0$, and 
	\begin{align}
		\| \Omega \|_{C^r}\leq C_r \|R\|_{C^{r+d+2}}.
	\end{align}
\end{Pro}

This result can be proved by using the preceding idea of smooth dependence on parameters. Indeed, thanks to Lemma \ref{Lem_base_tame2}, by slight adaptation of the arguments in subsection \ref{subsect_thetwistcase} we are able to obtain analogues of Lemma \ref{Lem_Mcont} and Lemma \ref{Lem_xieta} for the (untwisted) operator  $\cL_2$. Then Proposition \ref{Pro_LRS1111} follows by using arguments similar to the proof of Proposition \ref{Pro_LRS0}. So we will not repeat it here.

 We end this section by remarking that for the commuting maps $\FF=\TA+\bff$ and $\GG=\TB+\bfg$ in Section \ref{Section_actionofproducttype}, $\cL_1(\bff_1, \bfg_1)\neq 0$ and $\cL_2(\bff_2, \bfg_2)\neq 0$ in general (instead, they are quadratic, see Lemma \ref{Lem_comm_est}). So Proposition \ref{Pro_LRS0} and Proposition \ref{Pro_LRS1111} can not be applied directly to solve the cohomological equations in \eqref{twis_LE} and \eqref{untwis_LE}. Somehow, we will attempt to split $\bff_i, \bfg_i$, $i=1,2$, into $\bff_i=\cP(\bff_i)+\calE(\bff_i)$ and $\bfg_i=\cP(\bfg_i)+\calE(\bfg_i)$ in a tame way, such that $\cL_i(\cP(\bff_i), \cP(\bfg_i))=0$ and  $\calE(\bff_i), \calE(\bfg_i)$ are suitably small. This will be done in Section \ref{Section_tamesplit}.

\section{Tame Splitting}\label{Section_tamesplit}
The goal of this section is to prove Proposition \ref{Pro_split}. It implies that the perturbation can be split into two terms due to the commutation relations: one for which the linearized equations are solvable, and the other ``quadratically small'' with tame estimates. Different from Section \ref{Section_Smoothdependparamt}, our arguments here are based on a specific and explicit construction.  

\subsection{Construction}  

\begin{Lem}\label{lem_defproj}
	Let $A$ be an ergodic automorphism of $\TT^d$  and $\eta$ be a nonzero number in $\CC$.  For any function $f(x,y)\in C^\infty(\TT^d\times\TT^s,\CC)$, we define a function $\omega=\omega(f)$ as follows:
$
	\omega(x,y)\overset{\textup{def}}=\sum\limits_{n\in \ZZ^d} a_n(y)\, \expnx,
$
where for nonzero $n\neq 0$, 
\begin{equation}\label{anydef}   
	a_n(y)\overset{\textup{def}}=
	\begin{cases}
	     -\sum\limits_{l\geq 0}\eta^{-(l+1)}\, \hat f_{(A^*)^l\,n}(y), &  \textup{~if~} n\hookrightarrow E^u(A^*),\, E^c(A^*)\\ \\
		\quad\sum\limits_{l\leq-1}\eta^{-(l+1)}\, \hat f_{(A^*)^l\,n}(y), & \textup{~if~} n\hookrightarrow E^s(A^*)
	\end{cases}
\end{equation}
For $n=0$, $a_0(y)\overset{\textup{def}}=(\eta-1)^{-1}\, \hat f_0(y)$ if $\eta\neq 1$, and $a_0(y)\overset{\textup{def}}=0$ if $\eta=1$.
Then, $\omega$ is $C^\infty$ and 
\begin{align}\label{omgr}
	\|\omega\|_{C^r}\leq C_r \|f\|_{C^{r+d+2+\tau}},
\end{align}
where $\tau=(d+1)\frac{|\ln |\eta||}{\ln \rho}$ with $\rho>1$ the expansion rate for $A^*$, see \eqref{parhyp_split}.
\end{Lem}         
\begin{proof}
In order to prove that $\omega$ is $C^r$ differentiable  for each $r\in\NN$, we consider any multi-indices $\alpha\in\NN^d$ and $\beta\in\NN^s$ satisfying $|\alpha|+|\beta|=r$, and differentiate $\omega$ formally to get 
\begin{equation}\label{dklhr}
	\partial_x^\alpha\partial_y^\beta\omega=\sum_{n\in \ZZ^d} \partial_y^\beta a_n(y)\cdot (i2\pi)^{|\alpha|}\cdot n^\alpha\cdot \expnx.
\end{equation} 
Then we only need to show that the formal sum \eqref{dklhr} is absolutely convergent.

If $n\hookrightarrow E^u(A^*)$,  using $\partial_y^\beta\Big(\hat f_{(A^*)^l\,n} (y)\Big)=\hat{(\partial_y^\beta f)}_{(A^*)^l\,n}(y)$ and \eqref{parhyp_split} it follows   that
\begin{align*}
	|\partial_y^\beta a_n|\leq \sum\limits_{l\geq 0}|\eta|^{-(l+1)} \, \left|\hat{(\partial_y^\beta f)}_{(A^*)^l\,n}\right|\leq \sum\limits_{l\geq 0} \frac{|\eta|^{-(l+1)} \|\partial_y^\beta f\|_{C^k}}{\|(A^*)^l\,n\|^k}
	\leq & \sum\limits_{l\geq 0}  \frac{|\eta|^{-(l+1)}\|\partial_y^\beta f\|_{C^k}}{\|(A^*)^l\,\pi_u(n)\|^k}\\
	\leq & C \sum\limits_{l\geq 0}|\eta|^{-(l+1)} \, \frac{\|\partial_y^\beta f\|_{C^k}}{\rho^{kl}\|\pi_u(n)\|^k}\\
	\leq  & C' \frac{\|\partial_y^\beta f\|_{C^k}}{\|\pi_u(n)\|^k}\leq 3^k C'\frac{\|\partial_y^\beta f\|_{C^k}}{\|n\|^k}
\end{align*}
provided  $k>\frac{-\ln |\eta|}{\ln |\rho|}$. Indeed, the choice of $k$ ensures the convergence of $\sum_{l\geq 0} |\eta|^{-(l+1)}\rho^{-kl}$. 

Similarly, for  $n\hookrightarrow E^s(A^*)$, $|\partial_y^\beta a_n|
	\leq C_k\frac{\|\partial_y^\beta f\|_{C^k}}{\|n\|^k}$
	holds provided that  
	$k>\frac{\ln |\eta|}{\ln |\rho|}$.
	
When $n\hookrightarrow E^c(A^*)$, by the Katznelson lemma (see Lemma \ref{Lem_Katznelson}), $\|\pi_u(n)\|\geq \gamma \|n\|^{-d}$. Then,
\begin{align*}
	\|(A^*)^l n\|\geq\|(A^*)^l\,\pi_u(n)\|\geq C\rho^{l}\|\pi_u(n)\|\geq C\gamma  \rho^{l} \|n\|^{-d}\geq C\gamma  \rho^{l-l_0} \|n\|
\end{align*}
 for all  $l\geq l_0$ where $l_0=\left[\frac{(d+1)\ln\|n\|}{\ln \rho}\right]+1$. For $0\leq l\leq l_0-1$, by \eqref{parhyp_split} we get
\[
\|(A^*)^ln\|\geq \|(A^*)^l\pi_c(n)\|\geq C(1+l)^{-d} \|\pi_c(n)\|\ge\frac{C}{3}(1+l)^{-d} \|n\|.
\]
Thus we deduce that (only need to consider the worst case, $|\eta|<1$)
\begin{align*}
	|\partial_y^\beta a_n| 
		\leq & \sum_{l\geq 0} |\eta|^{-(l+1)} \, \left|\hat{(\partial_y^\beta f)}_{(A^*)^l\,n}\right|\\
		\leq & C'\|\partial^\beta_y f\|_{C^k}\left(\sum_{l=0}^{l_0-1}|\eta|^{-(l+1)}(1+l)^{dk}\|n\|^{-k}+ \sum_{l=l_0}^{+\infty}|\eta|^{-(l+1)}\rho^{-k(l-l_0)}\|n\|^{-k}\right)\\
			\leq	& C''\|\partial^\beta_y f\|_{C^k} \|n\|^{-k}\left(|\eta|^{-l_0}\cdot l_0^{dk+1}+ |\eta|^{-l_0}\sum_{i=0}^{+\infty}|\eta|^{-i}\rho^{-ki}\right)\\
		\leq	& C'''\|\partial^\beta_y f\|_{C^k} \|n\|^{-k}\cdot |\eta|^{-l_0}\left(l_0^{dk+1}+ c\right)\\
		\leq & C'''' \|\partial^\beta_y f\|_{C^k} \|n\|^{-k}\cdot  \|n\|^{\tau} 	\cdot(\ln\|n\|)^{dk+1}\leq  C_k \|\partial^\beta_y f\|_{C^k}\, \|n\|^{-k+\tau+1}
\end{align*}
provided that $k>\frac{-\ln |\eta|}{\ln \rho}$. Here, $\tau=(d+1)\frac{|\ln|\eta||}{\ln\rho}$. 

Therefore, by \eqref{dklhr} we can estimate that
\begin{align*}
	|\partial_x^\alpha\partial_y^\beta\omega|\leq (2\pi)^{|\alpha|}\sum_{n\in \ZZ^d} |\partial_y^\beta a_n| \cdot \|n\|^{|\alpha|}\leq   (2\pi)^{|\alpha|} \, C_k \sum_{n\in \ZZ^d} \| f\|_{C^{k+|\beta|}}\, \|n\|^{-k+\tau+1+|\alpha|}
\end{align*}
provided that $k>\frac{|\ln |\eta||}{\ln \rho}$. In particular, taking $k=|\alpha|+d+2+\tau$ we get 
\begin{align*}
	|\partial_x^\alpha\partial_y^\beta\omega|\leq C_r  \| f\|_{C^{|\alpha|+|\beta|+d+2+\tau}}.
\end{align*}
This holds for any multi-indices $\alpha\in\NN^d$ and $\beta\in\NN^s$ satisfying $|\alpha|+|\beta|=r$, so $\omega$ is $C^r$ for any $r\geq 0$, and satisfies estimate \eqref{omgr}.	
\end{proof}

Recall the following two operators:
\begin{equation*}
	\cL_1(f_1, g_1)=\Big(f_1(Bx,y)-Bf_1(x,y)\Big)-\Big(
	  g_1(Ax,y)-Ag_1(x,y)\Big)
\end{equation*}
for functions $f_1, g_1 : \TT^d\times\TT^s\to  \RR^d$. 
\begin{equation*}
	\cL_2(f_2, g_2)=\Big(f_2(Bx,y)-f_2(x,y)\Big)-\Big(
	 g_2(Ax,y)-g_2(x,y)\Big)
\end{equation*}
for functions $f_2, g_2 : \TT^d\times\TT^s\to  \RR^s$.

Then, the following tame splitting holds.   
\begin{Pro}[Tame Splitting]\label{Pro_split}
Suppose that $\cL_1(f_1, g_1)=\Phi_1$  and 
$\cL_2(f_2, g_2)=\Phi_2$, where  $ \Phi_1\in C^\infty(\TT^d\times\TT^s, \RR^d)$ and $\Phi_2\in C^\infty(\TT^d\times\TT^s, \RR^s)$. Then,
\begin{enumerate}[\bf(I)]
	\item  there 
	exists a splitting: 
	$
		f_1=\cP(f_1)+\calE(f_1)$ and $g_1=\cP(g_1)+\calE(g_1)$
	such that 
	\[\cL_1(\cP(f_1),\cP(g_1))=0,\qquad \cL_1(\calE(f_1),\calE(g_1))=\Phi_1\] 
	\begin{align*}
		\|\cP(f_1), \cP(g_1)\|_{C^r}\leq C_r \,\|f_1\|_{C^{r+\sigma_2}},\qquad \|\calE(f_1), \calE(g_1)\|_{C^r}\leq C_r \,\|\Phi_1\|_{C^{r+\sigma_2}}
	\end{align*}
		\item  there 
	exists a splitting: 
	$
	f_2=[f_2]+\cP(f_2)+\calE(f_2)$ and $g_2=[g_2]+\cP(g_2)+\calE(g_2)$ where
 $[f_2](y)$ and $[g_2](y)$ are the averages over the base $\TT^d$, and 
	\[\cL_2(\cP(f_2),\cP(g_2))=0,\qquad \cL_2(\calE(f_2),\calE(g_2))=\Phi_2\] 
	\begin{align*}
		\|\cP(f_2), \cP(g_2)\|_{C^r}\leq C_r \,\|f_2\|_{C^{r+\sigma_2}}, \qquad \|\calE(f_2), \calE(g_2)\|_{C^r}\leq C_r \,\|\Phi_2\|_{C^{r+\sigma_2}}
	\end{align*}
	Moreover, the averages over the base $[\cP(f_2)](y)=[\cP(g_2)](y)=[\calE(f_2)](y)=[\calE(g_2)](y)=0$.
\end{enumerate}
Here, the integer $\sigma_2$  depends only on the dimensions $d>0, s>0$ and  $A$ and $B$. 
\end{Pro}

\begin{proof}
\textbf{(I)} If $A$ and $B$ are semisimple, then by choosing a proper basis in which $A$ and $B$ are simultaneously diagonalize,  the system $\cL_1(f_1, g_1)=\Phi_1$ splits into finitely many one-dimensional equations of the following form
\begin{align}\label{bmwiuqhw}
	\big(\theta(Bx,y)-\mu \theta(x, y)\big)-\big(\psi(Ax,y)-\lambda \psi(x, y)\big)=\phi(x,y)
\end{align}
where $\lambda\neq 1, \mu\neq 1$ are a pair of eigenvalues of the ergodic $A$ and $B$, and $\theta, \psi, \phi$ belong to $C^\infty(\TT^d\times\TT^s,\RR)$.
For simplicity, we introduce the notation 
\[\Delta_A^\lambda\psi:=\psi(Ax,y)-\lambda \psi(x,y),\qquad \Delta_B^\mu\theta:=\theta(Bx,y)-\mu \theta(x,y).\]

Applying Lemma \ref{lem_defproj} with the number $\eta=\lambda$ and the  function $f=\theta$, we can construct a $C^\infty$ function $\omega(x,y)$  satisfying     
\begin{align}\label{nxkanal}
	\|\omega\|_{C^r}\leq C_r \|\theta\|_{C^{r+r_0}}
\end{align}
where $r_0=d+2+\tau$ with $\tau=(d+1)\frac{|\ln |\lambda||}{\ln \rho}$ and $\rho>1$ is the expansion rate for $A^*$. Now we construct the projections  
\begin{equation*}
	\cP(\theta)\overset{\textup{def}}=\Delta_A^\lambda\omega=\omega(Ax,y)-\lambda \omega(x,y),\qquad  \cP(\psi)\overset{\textup{def}}=\Delta_B^\mu\omega=\omega(Bx,y)-\mu \omega(x,y).
\end{equation*}
As $A$ commutes with $B$,   $\Delta_B^\mu \cP(\theta)-\Delta_A^\lambda\cP(\psi)=0$. Next, we set 
\begin{equation*}
	\calE(\theta)\overset{\textup{def}}=\theta-\cP(\theta),\qquad  \calE(\psi)\overset{\textup{def}}=\psi-\cP(\psi),
\end{equation*}
so by \eqref{bmwiuqhw} it satisfies 
\begin{equation}\label{afqqq}
	\Delta_B^\mu \, \calE(\theta)-\Delta_A^\lambda\, \calE(\psi)=\phi.
\end{equation}
Note that  all functions $\cP(\theta)$,
 $\cP(\psi)$, $\calE(\theta)$ and $\calE(\psi)$ above are $C^\infty$ on $\TT^d\times\TT^s$.
 
\noindent\textbf{Estimates for  $\cP(\theta)$ and $\cP(\psi)$.}  Obviously, \eqref{nxkanal} implies that      
\begin{equation}
	\|\cP(\theta), \cP(\psi)\|_{C^r}\leq C_r \|\theta\|_{C^{r+r_0}},
\end{equation}
where we enlarge the constant $C_r$ if necessary. 

\noindent\textbf{Estimates for  $\calE(\theta)$ and $\calE(\psi)$.} We write  $\calE(\theta)$ in the following form of Fourier series expansion
\begin{equation}\label{fexpan_calE}
	\calE(\theta)(x,y)=\sum\limits_{n\in \ZZ^d} \hat{\calE(\theta)}_n(y)\cdot \expnx,
\end{equation}
then the fact that $\calE(\theta)=\theta-\cP(\theta)=\theta-\Delta_A^\lambda\omega$ together with Lemma \ref{lem_defproj} implies
 $\hat{\calE(\theta)}_0=0$,  and for nonzero $n\neq 0$,
\begin{equation}\label{wotyoiy}
	\hat{\calE(\theta)}_n(y)=
	\begin{cases}
		\lambda \sum\limits_{l\in \ZZ}\lambda^{-(l+1)}\cdot \hat \theta_{(A^*)^l\,n}(y), &  \textup{if~} n\hookrightarrow E^s(A^*) \textup{~and~}  A^*n\hookrightarrow E^u(A^*), E^c(A^*) \\
		\quad 0, & \textup{otherwise}
	\end{cases}
\end{equation}

$\bullet$ Then it requires us to estimate the following $\sum\nolimits^A \hat \theta_{n}(y)$ for every $n$ satisfying $n\hookrightarrow E^s(A^*)$ and $A^*n\hookrightarrow E^u(A^*), E^c(A^*)$,
\begin{equation}\label{dfgoga}
	\sum\nolimits^A \hat \theta_{n}(y):=\sum\limits_{l\in \ZZ}\lambda^{-(l+1)}\cdot \hat \theta_{(A^*)^l\,n}(y).
\end{equation}
The formal sum on the right-hand side of \eqref{dfgoga} is absolutely convergent because $n\in \ZZ^d\setminus\{0\}$ has non-trivial projections to the expanding subspace and contracting subspace of the ergodic $A^*$.  

$\bullet$ Let us now estimate the size of \eqref{dfgoga} with respect to $\|n\|$, for every $n$ satisfying $n\hookrightarrow E^s(A^*)$ and $A^*n\hookrightarrow E^u(A^*), E^c(A^*)$. By \eqref{bmwiuqhw} one has $\Delta_A^\lambda\psi= \Delta_B^\mu\theta-\phi$, so the obstructions for $ \Delta_B^\mu\theta-\phi$ with respect to $A$ vanish, namely, for any nonzero $n\in\ZZ^d$,
\begin{equation}\label{eq_nisia}
	\sum\nolimits^A \hat\theta_{B^*n}(y)-\mu\sum\nolimits^A \hat\theta_n(y)-\sum\nolimits^A \hat\phi_n(y)=0.
\end{equation}
All formal sums in \eqref{eq_nisia} are absolutely convergent (the proof is the same with \eqref{dfgoga}).  
Therefore, by iterating \eqref{eq_nisia}  backward and  forward  we formally obtain
\begin{equation}\label{doublesum}
\sum\nolimits^A \hat\theta_n(y)=-\sum\nolimits^B_+\sum\nolimits^A \hat\phi_n(y)=\sum\nolimits^B_-\sum\nolimits^A \hat\phi_n(y),
\end{equation} 
where the notation 
\begin{equation*}
	\sum\nolimits^B_{\substack{+\\(-)}}\sum\nolimits^A \hat\phi_n:=\sum_{\substack{(l,k)\in  H^+\\ ( (l,k)\in H^-)}}\lambda^{-(l+1)}\cdot \mu^{-(k+1)}\,\,\hat\phi_{(A^*)^l(B^*)^k\,n}
\end{equation*}
with the sets $H^+=\{(l,k): l\in \ZZ,  k\geq 0\}$ and $H^-=\{(l,k):  l\in \ZZ, k< 0\}$. Consequently,  estimating \eqref{dfgoga} is equivalent to estimating the double sum in \eqref{doublesum}.

Since $n\hookrightarrow E^s(A^*)$ with $A^*n\hookrightarrow E^u(A^*), E^c(A^*)$,  due to condition \ref{condHR} in Section \ref{Section_actionofproducttype} we have: either for all non-zero $(l,k)\in H^+$ or for all non-zero $(l,k)\in H^-$, the following  polynomial estimates hold 
\begin{equation}\label{pollowbound}
	\|(A^*)^l(B^*)^k\,n\|\geq \frac{C}{|(l,k)|^{2d}} \|n\|,
\end{equation}
where $|(l,k)|:=\max\{|l|, |k|\}$. See \cite[page 1837]{Damjanovic_Katok10}. 

Without loss of generality, we suppose that \eqref{pollowbound} holds on $H^+$. Then, in order to estimate the $C^r$ norm of $\sum\nolimits^A \hat \theta_{n}(y)$, by \eqref{doublesum} it is equivalent to  investigating the following sum 
\begin{align*}
	\-\sum_{(l,k)\in H^+} \lambda^{-(l+1)}\cdot \mu^{-(k+1)}\,\,\hat\phi_{(A^*)^l(B^*)^k\,n}(y)
\end{align*}
We will choose a suitable $M>0$ and split the above sum into two parts $S_{<M}(\phi)$ and $S_{\geq M}(\phi)$: one is the finite sum on $H^+_{<M}=\{(l,k)\in H^+,  |(l,k)|< M \} $ and the other is the infinite sum on $H^+_{\geq M}=\{(l,k)\in H^+,  |(l,k)|\geq M \}$. For $H^+_{<M}$ we use the polynomial estimates  \eqref{pollowbound}, and for $H^+_{\geq M}$ we use the exponential  estimates in Lemma \ref{lem_Abexp}. More precisely, by Lemma \ref{lem_Abexp} one has
\begin{align}\label{hfgnvmow}
	\|(A^*)^l(B^*)^k n\|\geq C e^{|(l,k)|\kappa_0}\|n\|^{-d}  \geq C e^{(|(l,k)|-M)\kappa_0}\|n\|
\end{align}
where we choose the integer $M=\big[\frac{d+1}{\kappa_0}\ln \|n\|\big]+1$.

To estimate $S_{\geq M}(\phi)$, we set $m_0:=\max\{|\lambda|, |\mu|, |\lambda|^{-1},|\mu|^{-1}\}$, using \eqref{hfgnvmow} it follows that for any integers $r\geq 0$ and  $p> a=\big[\frac{2(d+1)}{\kappa_0}\ln |m_0|\big]+1$,
\begin{equation}\label{nannie}
\begin{aligned}
	\left\|S_{\geq M}(\phi)\right\|_{C^r(\TT^s)}\leq & C'\sum_{H^+_{\geq M}}|\lambda|^{-(l+1)} |\mu|^{-(k+1)}\frac{\|\phi\|_{C^{r+p}}}{\|(A^*)^l(B^*)^k\,n\|^p}\\
	\leq &C'' \|\phi\|_{C^{r+p}}\|n\|^{-p}\sum_{H^+_{\geq M}} m_0^{2|(l,k)|} \,e^{-(|(l,k)|-M)\kappa_0 p}\\
	\leq &C'' \|\phi\|_{C^{r+p}}\|n\|^{-p}\,m_0^{2M}\sum_{H^+_{\geq M}} \left(m_0^2\, e^{-\kappa_0 p}\right)^{|(l,k)|-M}\\
	\leq &C''' \|\phi\|_{C^{r+p}}\|n\|^{-p+a}\sum_{H^+_{\geq M}} \left(m_0^2\, e^{-\kappa_0 p}\right)^{|(l,k)|-M}\leq C_{r,p} \|\phi\|_{C^{r+p}}\|n\|^{-p+a}
\end{aligned}
\end{equation}

To estimate $S_{< M}(\phi)$, we use \eqref{pollowbound} and it follows that for any integers $r\geq 0$ and $p>a$,
\begin{equation}\label{fjslw}
\begin{aligned}
	\left\|S_{< M}(\phi)\right\|_{C^r(\TT^s)}\leq &C'\sum_{H^+_{< M}}|\lambda|^{-(l+1)} |\mu|^{-(k+1)}\frac{\|\phi\|_{C^{r+p}}}{\|(A^*)^l(B^*)^k\,n\|^p}\\
	\leq &C'' \|\phi\|_{C^{r+p}}\|n\|^{-p}\sum_{H^+_{< M}} m_0^{2|(l,k)|}\, |(l,k)|^{2d p}\\
	\leq &C'' \|\phi\|_{C^{r+p}}\|n\|^{-p} M^{2}m_0^{2M} M^{2dp}\\
	\leq &C''' \|\phi\|_{C^{r+p}}\|n\|^{-p+a} M^{2+2dp}\leq C_{r,p} \|\phi\|_{C^{r+p}}\|n\|^{-p+a+1}
\end{aligned}
\end{equation}

Thus, \eqref{nannie} together with \eqref{fjslw} gives the $C^r$-estimate for $\sum\nolimits^A \hat \theta_{n}(y)$. Combined with \eqref{wotyoiy}, we obtain
\begin{equation}\label{auukjfh}
	\left\|\hat{\calE(\theta)}_n(y)\right\|_{C^r(\TT^s)}\leq  C_{r,p} \|\phi\|_{C^{r+p}}\|n\|^{-p+a+1}
\end{equation}
for any $p> a=\big[\frac{2(d+1)}{\kappa_0}\ln |m_0|\big]+1$.

$\bullet$ Now, for any multi-indices $\alpha\in \NN^d$ and $\beta\in \NN^s$ satisfying $|\alpha|+|\beta|=r'$, we deduce from \eqref{fexpan_calE} and \eqref{auukjfh} that
\begin{align*}
	\|\partial_x^\alpha\partial_y^\beta\calE(\theta)\|_{C^0}\leq \sum_{n\in \ZZ^d} (2\pi\|n\|)^{|\alpha|}
	\left\|\hat{\calE(\theta)}_n(y)\right\|_{C^{|\beta|}(\TT^s)}\leq \tilde{C}_{|\beta|,p}  \|\phi\|_{C^{|\beta|+p}}\sum_{n\in \ZZ^d}  \|n\|^{-p+a+1+|\alpha|}.
\end{align*}
In particular, taking $p=a+|\alpha|+d+2$ we get
$
	\|\partial_x^\alpha\partial_y^\beta\calE(\theta)\|_{C^0}\leq C_{r'} \|\phi\|_{C^{r'+a+d+2}}.
$
Since it holds for any multi-indices $\alpha\in \NN^d$ and $\beta\in \NN^s$, we have
 \begin{align}\label{nncbaf}
	\|\calE(\theta)\|_{C^{r'}}\leq C_{r'} \|\phi\|_{C^{r'+a+d+2}}
 \end{align}
for every $r'\geq 0$. 

It remains to estimate $\calE(\psi)$. By \eqref{afqqq}, $\calE(\psi)$ satisfies the equation $\Delta_A^\lambda\, \calE(\psi)=\Delta_B^\mu \, \calE(\theta)-\phi$, so we infer from Proposition \ref{Pro_LRS0} and \eqref{nncbaf} that 
\[
	\|\calE(\psi)\|_{C^{r'}}\leq K_{r'} \|\Delta_B^\mu \, \calE(\theta)-\phi\|_{C^{r'+\sigma_1}}\leq\tilde {K}_{r'} (\|\calE(\theta)\|_{C^{r'+\sigma_1}}+\|\phi\|_{C^{r'+\sigma_1}})\leq  C_{r', \sigma_1}
	\|\phi\|_{C^{r'+a+d+2+\sigma_1}}
\]
for any $r'\geq 0$.

Therefore, we have finished the proof in the case of semisimple $A$ and $B$. If $A$ and $B$ are not semisimple we may have Jordan blocks, then instead of the one-dimensional  equation \eqref{bmwiuqhw}, for each Jordan block we would get a system of equations. However, analogous to \cite[Lemma 4.5]{Damjanovic_Katok10},  this system of equations can be studied inductively in finitely many steps, starting from an equation of the form \eqref{bmwiuqhw}. We will not repeat the arguments here.

\noindent \textbf{(II)} It can be proved in the same fashion as part \textbf{(I)}. In fact, the proof will be simpler because  the untwisted operator $\cL_2(f_2, g_2)=\Big(f_2(Bx,y)-f_2(x,y)\Big)-\Big(g_2(Ax,y)-g_2(x,y)\Big)$, and hence there is no Jordan blocks. Applying Lemma \ref{lem_defproj} with $\eta=1$ and $f=f_2$, we can construct a $C^\infty$ function $\omega\in C^\infty(\TT^d\times\TT^s,\RR^s)$. Then, we define 
 \begin{equation*}
	\cP(f_2)\overset{\textup{def}}=\omega(Ax,y)-\omega(x,y),\qquad  \cP(g_2)\overset{\textup{def}}=\omega(Bx,y)-\omega(x,y)
\end{equation*} 
and   
\begin{equation*}
	\calE(f_2)\overset{\textup{def}}=f_2-[f_2]-\cP(f_2),\qquad  \calE(g_2)\overset{\textup{def}}=g_2-[g_2]-\cP(g_2)
\end{equation*}
The remaining proof is just similar to that of part \textbf{(I)}, so we will not repeat it here.
\end{proof}

\subsection{Concluding remark}

Let $\rho:\ZZ^2\to \textup{Diff}^\infty(\TT^d\times\TT^s)$ denote the $\ZZ^2$ action generated by $\rho(\mathbf{e}_1)=\cT_{A,0}=A\times id_{\TT^s}$ and $\rho(\mathbf{e}_2)=\cT_{B,0}=B\times id_{\TT^s}$. Let $\cV$ be the space of all smooth functions $h=(h_1,h_2): \TT^d\times\TT^s\to \RR^d\times\RR^s$ which satisfy $\int_{\TT^d} h_2(x,y)\, dx=0$. We define two smooth tame linear operators $\Delta: \cV\to \cV\times \cV$ and $\cL: \cV\times\cV\to\cV$ by
 \begin{equation}\label{opedelcL}
 \begin{aligned}
 	\Delta h:=(\Delta^{\mathbf{e_1}} h,~ \Delta^{\mathbf{e_2}} h),\qquad
 	\cL(f, g):= \Delta^{\mathbf{e_2}}f -\Delta^{\mathbf{e_1}}g
 \end{aligned}
 \end{equation}
where the linear operators $\Delta^{\mathbf{e_1}}$ and $\Delta^{\mathbf{e_2}}$ are defined as follows:  for each $h=(h_1, h_2)\in \cV$ 
\begin{align*}
	\Delta^{\mathbf{e_1}}h:=h\circ \rho(\mathbf{e}_1) -\rho_*(\mathbf{e}_1) h &=(h_1(Ax,y)-A h_1(x,y),~ h_2(Ax,y)- h_2(x,y) ),\\
	\Delta^{\mathbf{e_2}}h:=h\circ \rho(\mathbf{e}_2) -\rho_*(\mathbf{e}_2) h &=(h_1(Bx,y)-B h_1(x,y), ~h_2(Bx,y)- h_2(x,y) ).
\end{align*}
Since $A$ commutes with $B$,  the following sequence is exact, i.e., $\cL\circ \Delta =0$,
\begin{equation}\label{exactsequence}
	\cV \xrightarrow{\quad\Delta\quad} \cV\times\cV \xrightarrow{\quad\cL\quad} \cV 
\end{equation} 
Let $Y=\textup{Im} \cL$ be the image of  $\cL$. As a consequence of Proposition \ref{Pro_split}, the following result holds.  
\begin{Cor}
	The exact sequence \eqref{exactsequence} admits the following tame splitting: there exist two smooth tame operators $\tilde\Delta:\cV\times\cV\to \cV$ and $\tilde\cL: Y\to\cV\times\cV$ such that $\Delta\circ\tilde\Delta+\tilde\cL\circ\cL=id$ on the space $\cV\times\cV$. Here, $Y=\cL(\cV\times\cV)$ is a linear subspace of $\cV$. 
\end{Cor}  
\begin{proof}
For any two elements $f=(f_1, f_2)$ and $g=(g_1, g_2)$ in $\cV$, note that  $\int_{\TT^d}f_2(x,y)\,dx$ $=$ $\int_{\TT^d}g_2(x,y)\,dx=0$,  then we can apply Proposition \ref{Pro_split} to split $f$ and $g$ into 
  \[f=\cP(f)+\calE(f), \qquad g=\cP(g)+\calE(g),\]
   where $\cP(f)=(\cP(f_1), \cP(f_2))$ and $\calE(f)=(\calE(f_1), \calE(f_2))$, and $\cP(g)=(\cP(g_1),\cP(g_2))$ and $\calE(g)=(\calE(g_1),\calE(g_2))$.
   Moreover, $\cL_1(\cP(f_1), \cP(f_2))=0$ and $\cL_2(\cP(f_2), \cP(g_2))=0$.
   This, combined with Propositions \ref{Pro_LRS0} and \ref{Pro_LRS1111}, implies that there exists a unique solution $h\in\cV$ satisfies the equation  
   $\Delta h = (\cP(f), \cP(g))$, where $\Delta$ is defined in \eqref{opedelcL}. As a consequence, we can define the operator  $\tilde\Delta: \cV\times\cV\to \cV$ by $\tilde\Delta (f, g)=h$. Clearly, it is a tame operator. 
   
On the other hand, for any $\Phi\in Y$, by definition there exists a pair $(f, g)$ such that $\cL(f,g)=\Phi$, in other words, $\cL_1(f_1, g_1)=\Phi_1$ and $\cL_2(f_2,g_2)=\Phi_2$. Then, applying the construction in Proposition \ref{Pro_split} we define the operator $\tilde \cL:Y\to\cV$ by $\tilde\cL(\Phi)=(f-\cP(f), g-\cP(g))$.

To show that $\tilde\cL$ is well defined, suppose that there is another pair $(f', g')$ such that $\cL(f',g')=\Phi$, then we claim that $(f'-\cP(f'), g'-\cP(g'))=(f-\cP(f), g-\cP(g))$. Indeed, we can write $f'=f+u$ and $g'=g+v$. Because $\cL(f',g')=\cL(f,g)$, we obtain $\cL(u,v)=0$. Thus, according to the construction in the proof of Proposition \ref{Pro_split},  $\cP(u)=u$ and $\cP(v)=v$. By linearity, we get $f'-\cP(f')=f+u-\cP(f)-\cP(u)=f-\cP(f)$ and $g'-\cP(g')=g-\cP(g)$. Therefore, $\tilde\cL$ is well defined. Finally,  $\Delta\circ\tilde\Delta+\tilde\cL\circ \cL=id_{\cV\times\cV}$ follows immediately from the above definitions.  
\end{proof}

\section{The KAM scheme and proof of Theorem \ref{Element_Thm1}}\label{Section_KAMscheme}

\subsection{The inductive step}\label{subsection_induclem}
We first establish the inductive step of the KAM scheme.

Note that at each step  of the iterative process we need to deal with cohomological equations of the form \eqref{twis_LE}--\eqref{untwis_LE}. According to Proposition \ref{Pro_LRS0} and Proposition \ref{Pro_LRS1111}, a phenomenon of the loss of regularity could happen. To overcome the fixed loss of derivatives, we use the smoothing operators $\rS_N$ for functions of $\TT^{d+s}$, and solve approximately the following (truncated) system:\begin{equation}\label{trun_twis_LE}
\begin{aligned}
	\bfh_1\circ \TA-A\bfh_1=\rS_N\bff_1,\qquad
	\bfh_1\circ \TB-B\bfh_1=\rS_N\bfg_1.
\end{aligned}	
\end{equation}
and 
\begin{equation}\label{trun_untwis_LE}
\begin{aligned}
\bfh_2 \circ \TA- \bfh_2=\rS_N\bff_2,\qquad
	 \bfh_2\circ \TB-\bfh_2=\rS_N\bfg_2.\
\end{aligned}
\end{equation}

To obtain approximate solutions we apply Proposition \ref{Pro_split} to $\rS_N\bff$ and $\rS_N\bfg$, and get the splitting
\begin{equation*}
	\rS_N\bff_1=\cP(\rS_N\bff_1)+\calE(\rS_N\bff_1),\qquad \rS_N\bfg_1=\cP(\rS_N\bfg_1)+\calE(\rS_N\bfg_1)
\end{equation*} 
\begin{equation*}
	\rS_N\bff_2=[\rS_N(\bff_2)]+\cP(\rS_N\bff_2)+\calE(\rS_N\bff_2),\qquad \rS_N\bfg_2=[\rS_N(\bfg_2)]+\cP(\rS_N\bfg_2)+\calE(\rS_N\bfg_2)
\end{equation*}
so that 
\[\cL_1(\cP(\rS_N\bff_1), \cP(\rS_N\bfg_1))=0,\qquad\cL_2(\cP(\rS_N\bff_2), \cP(\rS_N\bfg_2))=0\]
 with the averaged (w.r.t $\TT^d$) terms $[\cP(\rS_N\bff_2)](y)=[\cP(\rS_N\bfg_2)](y)=0$, and  
 \[\cL_1(\calE(\rS_N\bff_1), \calE(\rS_N\bfg_1))=\cL_1(\rS_N\bff_1, \rS_N\bfg_1),\qquad  \cL_2(\calE(\rS_N\bff_2), \calE(\rS_N\bfg_2))=\cL_2(\rS_N\bff_2, \rS_N\bfg_2).\] 
 Then, by Proposition \ref{Pro_LRS0} and Proposition \ref{Pro_LRS1111} the system 
\begin{equation}\label{eq_psnfg1}
\begin{aligned}
	\bfh_1\circ\TA-A\bfh_1=\cP(\rS_N\bff_1),\qquad \bfh_1\circ\TB-B\bfh_1=\cP(\rS_N\bfg_1).
\end{aligned}	
\end{equation}
and the system 
\begin{equation}\label{eq_psnfg2}
\begin{aligned}
  \bfh_2\circ\TA-\bfh_2=\cP(\rS_N\bff_2),\qquad \bfh_2\circ\TB-\bfh_2=\cP(\rS_N\bfg_2).\
\end{aligned}
\end{equation}
have a  solution  $\bfh=(\bfh_1, \bfh_2)$, $\bfh_1\in C^\infty(\TT^d\times\TT^s,\RR^d)$ and $\bfh_2\in C_0^\infty(\TT^d\times\TT^s,\RR^s)$.

In what follows, we set $\sigma=\max\{\sigma_1,\sigma_2, d+2\}$, where $\sigma_1$ and $\sigma_2$ are given in Proposition \ref{Pro_LRS0} and Proposition \ref{Pro_split}.	

\begin{Pro}\label{Pro_iterate}
Consider the commuting maps $\FF=\TA+\bff$ and $\GG=\TB+\bfg$ satisfying condition \ref{condIP}, where $A, B$ satisfy condition \ref{condHR} given in Section \ref{Section_actionofproducttype}. Then, for $N>0$ there exists a $C^\infty$ function $\bfh=(\bfh_1,\bfh_2)$ solving the system \eqref{eq_psnfg1}--\eqref{eq_psnfg2}, and it satisfies    
\begin{equation}\label{hhhhnorm}
	\|\bfh\|_{C^r}\leq  C_{r',r,\sigma}\, N^{r-r'+2\sigma}\| \bff,~\bfg\|_{C^{r'}},\qquad \text{for~}r\geq r'\geq 0.
\end{equation}

For the map defined by $H=id+\bfh$, if $\|\bfh\|_{C^1}\leq \frac{1}{4}$, then $H$  has a smooth inverse, and we obtain new $C^\infty$ maps $\widetilde\FF=H^{-1} \circ\FF \circ H$ and $\widetilde\GG=H^{-1}\circ\GG \circ H$. In addition, the new errors 
$ \tilde \bff=\widetilde\FF-\TA$ and $\tilde\bfg=\tilde\GG-\TB$ satisfy
\begin{align}
	\|\widetilde\bff, ~\widetilde\bfg \|_{C^0}  \leq & C_{r,\sigma}   \left( N^{2\sigma}\|\bff, \bfg\|_{C^1}\|\bff, \bfg\|_{C^0}+\frac{\|\bff, \bfg\|_{C^{\sigma+r+1}}^2}{N^r}+\frac{\|\bff, \bfg\|_{C^{\sigma+r}}}{N^r}\right), \quad \textup{for~}   r\geq 0,\label{wtf0norm}\\
	\|\widetilde \bff,~\widetilde \bfg\|_{C^r}\leq & C_{r,\sigma}\,\Big(1+N^{2\sigma}\| \bff,~\bfg\|_{C^r}\Big),\qquad \textup{for~}r> 0.\label{wtf_rnorm}
\end{align}
\end{Pro}

In the following proof, for simplicity we write  $\|u\|_{C^r}\ll \|v\|_{C^s}$ if there exists a constant $C>0$ independent of $u, v$ such that $\|u\|_{C^r}\leq C\|v\|_{C^s}$. We write $\|u\|_{C^r}\ll_{r,s} \|v\|_{C^s}$ when we want to stress that the constant $C$ depend on $r$ and $s$. 
\begin{proof}
By  Proposition \ref{Pro_LRS0}, Proposition \ref{Pro_LRS1111} and the analysis above, there exists a smooth $\bfh=(\bfh_1,\bfh_2)$ such that $\bfh_1\in C^\infty(\TT^{d}\times\TT^s,\RR^d)$ and $\bfh_2\in C_0^\infty(\TT^{d}\times\TT^s,\RR^s)$  are solutions to  \eqref{eq_psnfg1} and \eqref{eq_psnfg2} respectively. Moreover, combined with Proposition \ref{Pro_split} it follows that 
\begin{align}\label{bfhcr}
	 \|\bfh\|_{C^r}\ll_{r,\sigma} \|\cP(\rS_N\bff_1), \cP(\rS_N\bff_2) \|_{C^{r+\sigma}}\ll_{r,\sigma}  \|\rS_N\bff_1, \rS_N\bff_2 \|_{C^{r+2\sigma}} 
	\ll_{r,r',\sigma} N^{2\sigma+r-r'} \|\bff\|_{C^{r'}},
\end{align}
where for the last inequality we used Lemma \ref{Lem_trun}. This proves the desired estimate \eqref{hhhhnorm}.

Now, we form $H=id+\bfh:$ $(x,y)\mapsto (x+\bfh_1, y+\bfh_2)$. If $\|\bfh\|_{C^1}\leq \frac{1}{4}$,  by Lemma \ref{Apdix_pro1}  it is a smooth diffeomorphism. Under this conjugacy the original maps $\FF, \GG$ become 
\begin{align*}
	\tilde\FF\overset{\textup{def}}= H^{-1}\circ \FF\circ H,\qquad  \tilde\GG\overset{\textup{def}}= H^{-1}\circ \GG \circ H.
\end{align*}
In the sequel we will estimate the new errors $\tilde\bff=\tilde\FF-\TA$ and $\tilde\bfg=\tilde\GG-\TB$.

\noindent\textit{\large (I) Estimate of $\tilde\bff, \tilde\bfg$ in $C^0$.}
We first  estimate $\tilde\bff=(\tilde\bff_1, \tilde\bff_2)$.  By the relation $H\circ\tilde\FF$ $=$ $\FF\circ H$ it follows that $\tilde\bff_1=A \bfh_1+\bff_1\circ H-\bfh_1\circ\tilde\FF$ and  $\tilde\bff_2=\bfh_2+\bff_2\circ H-\bfh_2\circ\tilde\FF$. As $\bfh=(\bfh_1,\bfh_2)$ solves the system \eqref{eq_psnfg1}--\eqref{eq_psnfg2}, we thus obtain
\begin{equation}\label{tfexpre}
\begin{aligned}
	\tilde \bff_i=-\cP(\rS_N\bff_i)+\bff_i\circ H-(\bfh_i\circ\tilde\FF-\bfh_i\circ\TA),\qquad i=1,2
\end{aligned}
\end{equation}
Observe that $\bff=(\bff_1, \bff_2)$ has the following splitting  
\begin{align*}
\bff_1=&\rS_N\bff_1+\rR_N\bff_1=\cP(\rS_N\bff_1)+\calE(\rS_N\bff_1)+\rR_N\bff_1\\	
\bff_2=&\rS_N\bff_2+\rR_N\bff_2 =[\rS_N \bff_2]+\cP(\rS_N\bff_2)+\calE(\rS_N\bff_2)+\rR_N\bff_2
\end{align*}
 where $\rS_N$ and $\rR_N$ are the operators given in Lemma \ref{Lem_trun}. Hence, \eqref{tfexpre} becomes 
\begin{align}\label{wwleiw}
\begin{array}{lll}
	\tilde{\bff}_1= & &\calE(\rS_N\bff_1)+ \rR_N\bff_1+(\bff_1\circ H-\bff_1)-(\bfh_1\circ\tilde\FF-\bfh_1\circ\TA),\\\\
	\tilde{\bff}_2= &[\bff_2]-[\rR_N\bff_2]+&\calE(\rS_N\bff_2)+\rR_N\bff_2+(\bff_2\circ H-\bff_2)-(\bfh_2\circ\tilde\FF-\bfh_2\circ\TA).
\end{array}
\end{align}
During estimating  $\tilde\bff_2$,  the hard part is 
the averaged term $[\bff_2](y)=\int_{\TT^d}\bff_2(x,y)\,dx$ which is only of  order one  without further information. It is here that the intersection property comes into play,  causing this term to be of higher order. More precisely, as $\tilde\FF(x,y)=(Ax+\widetilde\bff_1, y+\widetilde\bff_2)$ satisfies the intersection property \ref{condIP}, it follows that for each point $y$,
\[\big(\TT^d\times\{y\}\big)  ~\cap ~\tilde\FF\big(\TT^d\times\{y\}\big) \neq \emptyset, \]
This implies that for every $y$, the map $x\mapsto \widetilde\bff_2(x,y)$ has  zeros, and hence we infer that
\begin{equation*}
	\|\tilde\bff_2\|_{C^0}\leq  2\| \tilde\bff_2-[\bff_2]\|_{C^0}.
\end{equation*}
This combined with \eqref{wwleiw} gives
 \begin{equation}\label{tildebff_2}
 \begin{split}
 	 	\|\tilde\bff_2\|_{C^0}\leq & 2\| \tilde\bff_2-[\bff_2]\|_{C^0}\\
	\leq & 2\Big( \|\calE(\rS_N\bff_2)\|_{C^0}+2\|\rR_N\bff_2\|_{C^0}+\|\bff_2\circ H-\bff_2\|_{C^0}+\|\bfh_2\circ \tilde\FF-\bfh_2\circ \TA\|_{C^0}\Big)\\
	\leq & 2\Big( \|\calE(\rS_N\bff_2)\|_{C^0}+2\|\rR_N\bff_2\|_{C^0}+\|\bff\|_{C^1}\|\bfh\|_{C^0}+\|\bfh\|_{C^1}\|\tilde \bff\|_{C^0}\Big).
 \end{split}
 \end{equation}
Meanwhile, \eqref{wwleiw} also gives the following preliminary estimate for $\tilde\bff_1$,
 \begin{align}\label{tildebff1}
	\|\tilde \bff_1\|_{C^0}\leq   \|\calE(\rS_N\bff_1)\|_{C^0}+\|\rR_N\bff_1\|_{C^0}+\|\bff\|_{C^1}\|\bfh\|_{C^0}+\|\bfh\|_{C^1}\|\tilde \bff\|_{C^0}.
\end{align}
As $\|\tilde\bff\|_{C^0}=\|\tilde\bff_1,~\tilde\bff_2\|_{C^0}$, \eqref{tildebff_2} and \eqref{tildebff1} together imply that
\begin{equation*}
		 	\|\tilde\bff\|_{C^0}
	\leq  2\Big( \|\calE(\rS_N\bff_1),~\calE(\rS_N\bff_2)\|_{C^0}+2\|\rR_N\bff\|_{C^0}+\|\bff\|_{C^1}\|\bfh\|_{C^0}+\|\bfh\|_{C^1}\|\tilde \bff\|_{C^0}\Big),
\end{equation*}
which yields
\begin{equation*}
		 	(1-2\|\bfh\|_{C^1})\cdot\|\tilde\bff\|_{C^0}
	\leq  2\Big( \|\calE(\rS_N\bff_1),~\calE(\rS_N\bff_2)\|_{C^0}+2\|\rR_N\bff\|_{C^0}+\|\bff\|_{C^1}\|\bfh\|_{C^0}\Big).
\end{equation*}
Using $\|\bfh\|_{C^1}\leq \frac{1}{4}$ yields
\begin{equation}\label{dawobff}
		 	\|\tilde\bff\|_{C^0}
	\leq  4\Big( \|\calE(\rS_N\bff_1),~\calE(\rS_N\bff_2)\|_{C^0}+2\|\rR_N\bff\|_{C^0}+\|\bff\|_{C^1}\|\bfh\|_{C^0}\Big).
\end{equation}
Let us estimate the three terms on the right-hand side of \eqref{dawobff}.
\begin{itemize}
	\item To estimate $ \|\calE(\rS_N\bff_1),~\calE(\rS_N\bff_2)\|_{C^0}$ we apply Proposition \ref{Pro_split} to $\rS_N\bff_1$, $\rS_N\bff_2$, and obtain
	\begin{align}\label{qo_sjfk}
		 \|\calE(\rS_N\bff_1)\|_{C^0}\ll_{\sigma}   \, \|\cL_1( \rS_N\bff_1, \rS_N\bfg_1  )\|_{C^\sigma}\, ,\qquad  \|\calE(\rS_N\bff_2)\|_{C^0}	\ll_{\sigma}  \|\cL_2( \rS_N\bff_2, \rS_N\bfg_2 )\|_{C^\sigma}
	\end{align}
	Note that
\begin{align*}
	\cL_1(\rS_N \bff_1, \rS_N\bfg_1)
	=  & \Delta^B(\bff_1-\rR_N\bff_1)-\Delta^A(\bfg_1-\rR_N\bfg_1)\\
	=&\cL_1(\bff_1, \bfg_1)-\Delta^B\rR_N\bff_1 + \Delta^A\rR_N\bfg_1\\
	=&\rS_N\Big(\cL_1(\bff_1, \bfg_1)\Big)+\rR_N \Big(\cL_1(\bff_1, \bfg_1)\Big)-\Delta^B\rR_N\bff_1 + \Delta^A\rR_N\bfg_1
\end{align*}
Then, invoking Lemma \ref{Lem_comm_est} and Lemma \ref{Lem_trun} it follows that, for any $r\geq 0$,
\begin{align}
	\|\cL_1(\rS_N \bff_1, \rS_N\bfg_1)\|_{C^\sigma}\ll_{\sigma, r}\, & N^\sigma\|\cL_1(\bff_1, \bfg_1)\|_{C^0}+\frac{\|\cL_1(\bff_1, \bfg_1)\|_{C^{\sigma+r}}}{N^r}+\frac{\|\bff_1, \bfg_1\|_{C^{\sigma+r}}}{N^r}\nonumber\\
	 \ll_{\sigma, r}\, & N^\sigma\|\bff, \bfg\|_{C^1}\|\bff, \bfg\|_{C^0}+\frac{\|\bff, \bfg\|_{C^{\sigma+r+1}}^2}{N^r}+\frac{\|\bff, \bfg\|_{C^{\sigma+r}}}{N^r}.\nonumber
\end{align}
This together  with \eqref{qo_sjfk} gives
 \begin{align}\label{caleNf1}
	 \|\calE(\rS_N\bff_1)\|_{C^0}\ll_{\sigma,r}  \, N^\sigma\|\bff, \bfg\|_{C^1}\|\bff, \bfg\|_{C^0}+\frac{\|\bff, \bfg\|_{C^{\sigma+r+1}}^2}{N^r}+\frac{\|\bff, \bfg\|_{C^{\sigma+r}}}{N^r}.
\end{align}
	Similarly, we can also show that
\begin{align}\label{caleNf2}
	 \|\calE(\rS_N\bff_2)\|_{C^0}	\ll_{\sigma, r}\,  	  N^\sigma\|\bff, \bfg\|_{C^1}\|\bff, \bfg\|_{C^0}+\frac{\|\bff, \bfg\|_{C^{\sigma+r+1}}^2}{N^r}+\frac{\|\bff, \bfg\|_{C^{\sigma+r}}}{N^r}.
\end{align}
\item Using  Lemma \ref{Lem_trun},
\begin{equation}\label{es_RNf1f2}
	\|\rR_N\bff\|_{C^0}\ll_{r} \|\bff\|_{C^r}N^{-r}.
\end{equation}
\item Applying \eqref{bfhcr} with $r= r'=0$, we obtain
\begin{equation}\label{es_f1h0}
	\|\bff\|_{C^1}\|\bfh\|_{C^0}\ll_{\sigma}\, N^{2\sigma}\|\bff\|_{C^1}\|\bff\|_{C^0}.
\end{equation}
\end{itemize}
Therefore,  combining \eqref{caleNf1}--\eqref{es_f1h0} we have the following estimate for $\|\tilde\bff\|_{C^0}$,
\begin{equation}
		 	\|\tilde\bff\|_{C^0}
	\ll_{\sigma,r} N^{2\sigma}\|\bff, \bfg\|_{C^1}\|\bff, \bfg\|_{C^0}+\frac{\|\bff, \bfg\|_{C^{\sigma+r+1}}^2}{N^r}+\frac{\|\bff, \bfg\|_{C^{\sigma+r}}}{N^r}.
\end{equation}

On the other hand, to estimate $\tilde\bfg$ we repeat similar arguments as above and obtain
\begin{equation}
		 	\|\tilde\bfg\|_{C^0}
	\ll_{\sigma,r} N^{2\sigma}\|\bff, \bfg\|_{C^1}\|\bff, \bfg\|_{C^0}+\frac{\|\bff, \bfg\|_{C^{\sigma+r+1}}^2}{N^r}+\frac{\|\bff, \bfg\|_{C^{\sigma+r}}}{N^r}.
\end{equation}
This proves the desired estimate \eqref{wtf0norm}.

\noindent\textit{\large (II) Estimate of $\tilde\bff, \tilde\bfg$ in $C^r$.}
Note that $\tilde\bff$ can be expressed as follows 
\begin{align*}
	\tilde\bff=H^{-1}\circ \FF\circ H-\TA
	=& (H^{-1}-id)\circ\FF\circ H +(A\bfh_1, \bfh_2)+\bff\circ H,
\end{align*}
which leads to
\begin{equation*}
	\|\tilde\bff\|_{C^r}\ll \|(H^{-1}-id)\circ\FF\circ H\|_{C^r} +\|\bfh\|_{C^r}+\|\bff\circ H\|_{C^r}.
\end{equation*}
Note that $\tilde\bff(x,y)$ is a $\ZZ^{d+s}$-periodic function. By Lemma \ref{Apdix_linecont} it follows that
\begin{align*}
	\|(H^{-1}-id)\circ\FF\circ H\|_{C^r}\ll_{r} \, 1+\|(H^{-1}-id)\|_{C^r}+\|\bff\|_{C^r}+\|\bfh\|_{C^r}
\end{align*}
and
\begin{align*}
	\|\bff\circ H\|_{C^r}\ll_r\, 1+\|\bff\|_{C^r}+\|\bfh\|_{C^r}
\end{align*}
By Lemma \ref{Apdix_pro1}, $\|H^{-1}-id\|_{C^r}\ll_r \|\bfh\|_{C^r}$. 
Thus, combined with \eqref{bfhcr} we obtain
\begin{align*}
	\|\tilde\bff\|_{C^r}\ll_r\, 1 +\|\bfh\|_{C^r}+\|\bff\|_{C^r}\ll_r\, 1+N^{2\sigma} \|\bff\|_{C^r}.
\end{align*}

Similarly, $\|\tilde\bfg\|_{C^r}$ can be estimated in the same way as $\|\tilde\bff\|_{C^r}$. This finally finishes the proof.
	
\end{proof}

\subsection{Proof of Theorem \ref{Element_Thm1}}\label{subsection_KAMscheme}

Based on Proposition \ref{Pro_iterate}, we now prove Theorem \ref{Element_Thm1}.

\begin{proof}[Proof of Theorem \ref{Element_Thm1}]
To begin the iterative process, we set up
\[\FF^{(0)}=\TA+\bff^{(0)},\quad \GG^{(0)}=\TB+\bfg^{(0)}, \quad H^{(0)}=id\]
where $\bff^{(0)}=\bff$ and $\bfg^{(0)}=\bfg$.
Fix  a sufficiently large integer $N_0>0$, and define $N_i$ inductively by
\begin{align}\label{Nisetting}
	N_{i+1}=N_{i}^{\frac{3}{2}}
\end{align}
for all $i=0,1,2,\cdots$.

We will construct $ \FF^{(i)},\GG^{(i)}, H^{(i)}$ inductively for $i\geq 0$: suppose that we have already obtained $\FF^{(i)},\GG^{(i)}, H^{(i)}$, then at  the $(i+1)$-th step we choose the smoothing operator $\rS_{N_i}$ with  $N_{i}>0$ given in \eqref{Nisetting}  and apply  Proposition \ref{Pro_iterate} to obtain $\bfh^{(i)}$.
This produces a new smooth conjugacy
\[H^{(i+1)}=id+\bfh^{(i)}.\]
Let us use the notation
\begin{align*}
\cE_{i,r}\overset{\textup{def}}=\|\bff^{(i)},~\bfg^{(i)}\|_{C^r}\,, \qquad
 \cU_{i,r}\overset{\textup{def}}=\|\bfh^{(i)}\|_{C^r}.  
\end{align*}
Then by Proposition \ref{Pro_iterate} we have
\begin{align}
	\cU_{i,r}\ll_{r,r',\sigma}  \, N_{i}^{r-r'+2\sigma}\,\cE_{i,r'},\qquad \text{for~}r\geq r'\geq  0.\label{iterate_h_rnorm}
\end{align}

To proceed, we have to check that each $H^{(i+1)}$ is indeed invertible. This is true if  the following condition is satisfied 
\begin{equation}\label{smalonc1}\tag{D}
	\cU_{i,1}\leq\frac{1}{4}.
\end{equation}
Then we can define
\[\FF^{(i+1)}=\left(H^{(i+1)}\right)^{-1}\circ \FF^{(i)}\circ H^{(i+1)}=\TA+\bff^{(i+1)},\]
 \[\GG^{(i+1)}=\left(H^{(i+1)}\right)^{-1}\circ \GG^{(i)}\circ H^{(i+1)}=\TB+\bfg^{(i+1)},\]  
 where $\bff^{(i+1)}$ and $\bfg^{(i+1)}$ are the new errors. 
 According to Proposition \ref{Pro_iterate},  
\begin{align}
	\cE_{i+1,0} \ll_{r,\sigma} &\, N_{i}^{2\sigma}\,\cE_{i,1}\cdot\cE_{i,0}+\frac{\cE^2_{i,\sigma+r+1}}{ N_{i}^r}+\frac{\cE_{i,\sigma+r}}{ N_{i}^r},\qquad \text{for~} r\geq 0.\label{iterate_wtf0norm}\\
\cE_{i+1,r}\ll_{r,\sigma}  &\, 1+ N_{i}^{2\sigma}\,\cE_{i,r},\qquad \text{for~}r> 0.\label{iterate_wtf_rnorm}
\end{align}

To ensure the above condition \eqref{smalonc1} and the convergence of the KAM scheme, one only needs that for the original error, $\cE_{0,0}=\|\bff^{(0)},~\bfg^{(0)}\|_{C^0}$ is sufficiently small and  $\cE_{0,\mu_0}=\|\bff^{(0)},~\bfg^{(0)}\|_{C^{\mu_0}}$ is well controlled for some integer $\mu_0>0$. This is guaranteed by the 
 following lemma.
\begin{Lem}\label{Lem_induc_ineq}
Let us set $\mu_0=20(\sigma+1)$ and
$\kappa=6(\sigma+1)$. We can choose  $N_0>0$ suitably large, such that if
$\cE_{0,0}=\|\bff^{(0)},~\bfg^{(0)}\|_{C^0}\leq N_0^{-\kappa}$ and
 $\cE_{0,\mu_0}=\|\bff^{(0)},~\bfg^{(0)}\|_{C^{\mu_0}}\leq N_0^{\frac{3}{4}\kappa}$, then  condition \eqref{smalonc1} holds for all $i\geq 0$. In addition, we have
	\begin{align}\label{3fchi}
 \cE_{i,0}\leq N_i^{-\kappa},\qquad
 \cE_{i,\mu_0}\leq N_i^{\frac{3}{4}\kappa},\qquad \cU_{i,1}\leq N_i^{-\frac{1}{2}\kappa}
\end{align}
\end{Lem}
\begin{proof}
We prove it by induction. For $i=0$, by \eqref{iterate_h_rnorm} we have $\cU_{0,1}\ll N_0^{2\sigma+1}\cE_{0,0}$ $\leq N_0^{2\sigma+1-\kappa}$. As $N_0>0$ is large enough, it follows that $\cU_{0,1}\leq N_0^{-\frac{\kappa}{2}}$. Thus, by assumption, \eqref{3fchi} holds for $i=0$, and condition  \eqref{smalonc1} for $i=0$ is satisfied provided that  $N_0$ is large.

Suppose that \eqref{3fchi} holds for all steps $\leq i$, we need to verify these estimates for the $(i+1)$-th step. Note that by the interpolation inequalities (Lemma \ref{cor_intpest}), we get
\begin{equation}\label{intp_inequ}
	\cE_{i,1}\leq C_{\mu_0}\,\cE_{i,0}^{1-\frac{1}{\mu_0}} \,\cE_{i,\mu_0}^{\frac{1}{\mu_0}}\leq C_{\mu_0}\,N_i^{-\kappa(1-\frac{7}{4\mu_0})}.
\end{equation}
Then, 
using inequality \eqref{iterate_wtf0norm} with $r=\mu_0-\sigma-1$ we  obtain 
\begin{align}
	\cE_{i+1,0}\ll_{\sigma}  N_{i}^{2\sigma}\cE_{i,1}\cdot\cE_{i,0}+\frac{\cE^2_{i,\mu_0}}{N_{i}^{\mu_0-\sigma-1}}+\frac{\cE_{i,\mu_0-1}}{N_{i}^{\mu_0-\sigma-1}} & \ll_{\sigma}   \, N_{i}^{2\sigma}N_i^{-\kappa(2-\frac{7}{4\mu_0})} +\frac{N_i^{\frac{3}{2}\kappa}}{N_{i}^{\mu_0-\sigma-1}} \nonumber\\
	 & \ll_{\sigma}   N_{i}^{2\sigma-\kappa(2-\frac{7}{4\mu_0})} + N_{i}^{-10(\sigma+1)}\nonumber\\
	 &< \, N_{i+1}^{-\kappa}\,\label{cE_i1}
\end{align}
provided that  $N_0$ is suitably large. Next, applying \eqref{iterate_wtf_rnorm} with  $r=\mu_0$, it follows that
\begin{align}\label{rho}
	\cE_{i+1,\mu_0}\ll_{\mu_0,\sigma} 1+N_{i}^{2\sigma}\cE_{i,\mu_0} \ll_{\mu_0,\sigma} N_{i}^{2\sigma} N_i^{\frac{3}{4}\kappa}
	< N_{i+1}^{\frac{3}{4}\kappa}.
\end{align}
Finally, applying inequality  \eqref{iterate_h_rnorm}  with $r=1$ and $r'=0$,  we have
\begin{align}\label{cu_i11}
	\cU_{i+1,1}	&\ll_{\sigma} N_{i+1}^{2\sigma+1}\cE_{i+1,0}  \ll_{\sigma}   N_{i+1}^{2\sigma+1}N_{i+1}^{-\kappa}<
	N_{i+1}^{-\frac{1}{2}\kappa}.
\end{align}
This also implies that condition \eqref{smalonc1} holds at the $(i+1)$-th step because $N_{i+1}=N_0^{\left(\frac{3}{2}\right)^{i+1}}$ is large.
This proves Lemma \ref{Lem_induc_ineq}.
\end{proof}

Now, let us proceed with the proof of Theorem \ref{Element_Thm1}. By Lemma \ref{Lem_induc_ineq}, as long as $\cE_{0,0}$ is suitably small and
 $\cE_{0,\mu_0}\leq 1$, the following sequences 
\[\|\bff^{(i)},~\bfg^{(i)}\|_{C^0}\leq N_0^{-\left(\frac{3}{2}\right)^i\kappa},\quad \|\bfh^{(i)}\|_{C^1}\leq  N_0^{-\frac{1}{2}\left(\frac{3}{2}\right)^i\kappa} \]
converge rapidly to zero. This rapid convergence ensures  that 
as $l\to\infty$, the composition 
\[\mathcal{H}_l =H^{(1)}\circ\cdots\circ H^{(l)}\] converges in the $C^1$ topology to  some $\mathcal{H}_\infty$ which is a $C^1$ diffeomorphism, for which the following equations hold
\[ \FF\circ \mathcal{H}_\infty=\mathcal{H}_\infty \circ \TA,\qquad \GG\circ \mathcal{H}_\infty=\mathcal{H}_\infty\circ \TB.\]

It remains to show that the above $C^1$ limit solution $\mathcal{H}_\infty$ is also of class $C^p$ for any $p>1$. In fact, as shown in \cite{Zeh_generalized1}, this can be achieved by making full use of the interpolation inequalities. More precisely, observe that for  any $m>0$, applying \eqref{iterate_wtf_rnorm} with $r=m$ we get
\[
	\cE_{i,m}\leq C_{m,\sigma} \big(1+N_{i-1}^{2\sigma}\,\cE_{i-1,m}\big)
\]
for some constant $C_{m,\sigma}>1$. This also gives $1+\cE_{i,m}\leq C_{m,\sigma} \, N_{i-1}^{2\sigma}\Big(1+\cE_{i-1,m}\Big)$, from which we derive inductively that
\begin{align}\label{cnwnkw}
	\cE_{i,m}\leq \left(1+\cE_{0,m} \right)\prod_{j=0}^{i-1}\left(C_{m,\sigma} \,N_{j}^{2\sigma}\right)
	\leq  \left(1+\cE_{0,m} \right) C^i_{m,\sigma}\left(\prod_{j=0}^{i-1} N_j\right)^{2\sigma}
	&\leq  M_m\cdot C^i_{m,\sigma}\cdot N_0^{\left(\frac{3}{2}\right)^i 4\sigma}\nonumber\\
	&= M_m\cdot C^i_{m,\sigma}\cdot N_i^{4\sigma}
\end{align}
where we denote $M_m= (1+\cE_{0,m})>1$.

In particular, for any given $p>1$, we choose $m=4p$ and apply the interpolation inequalities (Lemma \ref{cor_intpest}) to obtain that
\begin{align*}
	\cE_{i,p}\leq  C_{p}\,  \cE_{i,m}^{\frac{1}{4}}\cdot \cE_{i,0}^{\frac{3}{4}}\leq
	&C_p\left(M_m\, C^i_{m,\sigma}\right)^{\frac{1}{4}}\,   N_i^{\sigma}\cdot N_{i}^{-\frac{3}{4}\kappa}
	\leq   C_p\, M_m\, C_{m,\sigma}^{i}\cdot N_i^{-\frac{1}{2}\kappa}
\end{align*}
where we have used \eqref{3fchi} and \eqref{cnwnkw}. Then, applying  \eqref{iterate_h_rnorm} with $r=r'=p$ yields 
\begin{align*}
\cU_{i,p}\leq  C_{p,\sigma}\, N_{i}^{2\sigma}\cE_{i,p}
\leq  & C_{p,\sigma} C_p\,M_m\, C^{i}_{m,\sigma}\cdot N_i^{2\sigma-\frac{\kappa}{2}}
\leq L \cdot b^i\cdot N_i^{-1-\sigma}
\end{align*}
where the constants $L=C_{p,\sigma} C_p\,M_m$  and $b=C_{m,\sigma}>1$, with $m=4p$. Observe that although $b^i$ grows exponentially, the quantity $N_i^{-1-\sigma}$ decays super-exponentially. Hence,  $\cU_{i,p}=\|\bfh^{(i)}\|_{C^p}\ll N_i^{-\sigma}$ converges rapidly to zero, which ensures the convergence of the sequence $\{\mathcal{H}_l\}_l$ in the $C^p$ topology, and the limit is exactly $\mathcal{H}_\infty$.  Since the above argument is true for any given integer $p\geq 1$,  we conclude that $\mathcal{H}_\infty$ is $C^\infty$. This finally finishes the proof of Theorem \ref{Element_Thm1}.
\end{proof}

\section{Proofs of Theorems \ref{MainThm_0}, \ref{MainThm_1} and \ref{MainThm_2}}\label{Section_proofMainResult}

\subsection{Proof of Theorem \ref{MainThm_0}}\label{subsecproofA}

Based on Proposition \ref{Pro_conj_ave} and Theorem \ref{Element_Thm1} we proceed to prove Theorem \ref{MainThm_0}.

We will also need the following lemma.

\begin{Lem}\label{lem_periodid}
	Given $\theta\in \QQ^s$ and an integer $q>0$ satisfying $q\theta\in \ZZ^s$. Suppose that  $\Phi(y):\TT^s\to\TT^s$ is a $C^\infty$ diffeomorphism which satisfies:
	\begin{enumerate}
		\item $\Phi$ is  $C^1$-sufficiently close to the toral translation  $R_\theta: y\mapsto y+\theta$ (mod $\ZZ^s$). More precisely, we can write $\Phi=R_\theta+\omega$ where the function $\omega\in C^\infty(\TT^s,\RR^s)$ and $\|\omega\|_{C^1}\ll 1$.
		\item  the $q$-fold composition $\Phi^q$ of $\Phi$  satisfies that $\Phi^q=id_{\TT^s}$
	\end{enumerate}
Then, $\Phi$ can be $C^\infty$-conjugated to $R_\theta$ via a conjugacy $V=id_{\TT^s}+v$ with $v$ of the form 
\begin{equation*}
	v(y)=\frac{1}{q}\sum_{i=0}^{q-2}(q-i-1)\,\omega\circ\Phi^{i}(y),
\end{equation*}  
and hence $\|v\|_{C^1}\leq C\|\omega\|_{C^1}$. 
\end{Lem}
\begin{Rem}
	The conjugacy between $\Phi$ and $R_{\theta}$ is not unique, here we just construct one having an explicit formula. The non-uniqueness is caused by the non-ergodicity of $\Phi$. 
\end{Rem}

\begin{proof}
By assumption  the diffeomorphism  $\Phi$ can be written as $\Phi=R_\theta+\omega$, where $\omega\in C^\infty(\TT^s,\RR^s)$ with $\|\omega\|_{C^1}$ sufficiently small. Note that $\Phi^i(y)=\Phi\circ\Phi^{i-1}(y)$ $=$ $\Phi^{i-1}(y)+\theta+\omega\circ \Phi^{i-1}(y)$. By iterating this formula for $i=1,\cdots,q$ and adding them up, one gets 
\begin{equation*}
	\Phi^q(y)=y+q\theta+\sum_{i=0}^{q-1}\omega\circ\Phi^{i}(y)=R_{q\theta}+\sum_{i=0}^{q-1}\omega\circ\Phi^{i}.
\end{equation*}	
Due to $q\theta\in \ZZ^s$ and $\Phi^q=id_{\TT^s}$, we have $\sum_{i=0}^{q-1}\omega\circ\Phi^{i}=0$ mod $\ZZ^s$. As $\omega$ is sufficiently small, it implies that 
\begin{align}\label{qomegzero}
	\sum_{i=0}^{q-1}\omega\circ\Phi^{i}=0.
\end{align}

Now, we proceed to construct a near-identity conjugacy $V\in \textup{Diff}^\infty(\TT^s)$ such that  $V\circ \Phi=R_\theta\circ V$. We write $V(y)=y+v(y)$ with $v\in C^\infty(\TT^s,\RR^s)$, then the conjugacy equation $V\circ \Phi=R_\theta\circ V$ reduces to the following equation
\begin{equation}\label{coeqvomega}
	v(y)-v\circ\Phi(y)=\omega(y).
\end{equation}
Using \eqref{qomegzero}, it is easy to check that equation \eqref{coeqvomega} has a $C^\infty$ solution given by
\begin{equation}
	v(y):=\frac{1}{q}\sum_{i=0}^{q-2}(q-i-1)\,\omega\circ\Phi^{i}(y).
\end{equation}
The condition $\|\omega\|_{C^1}\ll 1$ ensures that $\|v\|_{C^1}$ is also suitably small, so the map $V=id_{\TT^s}+v$ has a smooth inverse. Therefore, $\Phi$ is $C^\infty$-conjugate to $R_\theta$ via the near-identity conjugacy $V$. 
\end{proof}

\begin{proof}[Proof of Theorem \ref{MainThm_0}]
Thanks to Proposition \ref{Pro_conj_ave}, the original perturbation $\alpha=\langle \cT_{A_1,\tau_1}, \cT_{A_2,\tau_2} \rangle$ is $C^\infty$-conjugate to the action $\langle \cT_{A_1,[\tau_1]}, \cT_{A_2,[\tau_2]} \rangle$. So we only need to prove Theorem \ref{MainThm_0} for the case where $\tau_1(x)$ and $\tau_2(x)$ are constants, i.e., $\tau_1=[\tau_1]$ and $\tau_2=[\tau_2]$, and for brevity we will denote  
\begin{equation*}
	\theta_1:=[\tau_1],\qquad \theta_2:=[\tau_2]
\end{equation*}
and from now on, our unperturbed action is assumed to be $\alpha=\langle\cT_{A_1,\theta_1}, \cT_{A_2,\theta_2}\rangle$, where $\theta_1, \theta_2\in \QQ^s$. In other words,  $\alpha(\mathbf{n})=\cT_{A_1^{n_1}A_2^{n_2}, n_1\theta_1+n_2\theta_2}$ for any $\mathbf{n}=(n_1,n_2)\in \ZZ^2$.

 \textsl{In the case of integer $(\theta_1, \theta_2)$.} If $(\theta_1, \theta_2)\in \ZZ^s\times\ZZ^s$, the unperturbed action $\alpha$ becomes $\alpha=\langle\cT_{A_1,0},\cT_{A_2,0}\rangle$, then our result follows immediately from Theorem \ref{Element_Thm1}.

 \textsl{In the case of non-integer $(\theta_1, \theta_2)$.} If $(\theta_1,\theta_2)\in \QQ^s\times\QQ^s$ with $(\theta_1,\theta_2)\notin \ZZ^s\times\ZZ^s$, we split the proof into two parts.
 
\noindent\textbf{\large Part 1.} Recall that  $M_0\geq 1$ denotes the minimal positive integer $\lambda$ such that $\lambda\,\theta_1\in \ZZ^s$, $\lambda\, \theta_2\in \ZZ^s.$ In the sequel we denote
\begin{align*}
	\ZZ^2_{M_0}:=\{(i,j)\in \ZZ^2: |i|\leq M_0,~ |j|\leq M_0 \},
\end{align*}
and consider two sets as follows
\begin{align*}
\Sigma:=\{ (i,j)\in \ZZ^2_{M_0}: i\theta_1+j\theta_2\in \ZZ^s\},\qquad\Lambda:=\{ (i,j)\in \ZZ^2_{M_0}: i\theta_1+j\theta_2\notin \ZZ^s\}.
\end{align*}
$\Sigma\neq \emptyset$ and $\Lambda\neq\emptyset$ since $(\theta_1,\theta_2)$ is non-integer. So there exists $\delta^*=\delta^*(\theta_1,\theta_2)>0$ such  that 
	\begin{equation}\label{mindistij}
		 \textup{dist}(i\theta_1+j\theta_2, \ZZ^s)\geq \delta^*,\qquad \forall (i,j)\in \Lambda
	\end{equation}
	since $\Lambda$ is a finite set. Now, we consider the perturbed action $\tilde\alpha$. For $r\geq 0$ we denote by
\begin{equation*}
	d_{{C^r}}(\tilde\alpha, \alpha; M_0):=\max_{\mathbf{k}\in \ZZ^2_{M_0}} \textup{dist}_{C^r}(\tilde\alpha(\mathbf{k}), \alpha(\mathbf{k})).
\end{equation*}
the maximum of $\textup{dist}_{C^r}(\tilde\alpha(\mathbf{k}), \alpha(\mathbf{k}))$ for all $\mathbf{k}\in \ZZ^2_{M_0}$. Then, the following properties hold:
		
 \textbf{(I)} \textit{If 
 \begin{equation}\label{swiqon}
d_{C^0}(\tilde\alpha, \alpha; M_0)<\frac{\delta^*}{2},
\end{equation}
 and if for some $\mathbf{k}\in\ZZ^2_{M_0}$ the map $\tilde\alpha(\mathbf{k})$ satisfies the intersection property \ref{condIP}, then $\mathbf{k}\in \Sigma$}.

We now argue by contradiction. Assume that $\mathbf{k}\in \Lambda$. Taking a $d$-dimensional torus $\Gamma=\TT^d\times\{y=0\}$,  by \eqref{mindistij} we see that the Hausdorff distance between $\Gamma$ and its image under the map $\alpha(\mathbf{k})=\cT_{A_1^{k_1}A_2^{k_2}, k_1\theta_1+k_2\theta_2}$ is greater than $\delta^*$. This, together with \eqref{swiqon}, implies that the Hausdorff distance between  $\Gamma$ and its image under  $\tilde\alpha(\mathbf{k})$ is greater than $\frac{\delta^*}{2}$, so they cannot intersect.  This contradicts the intersection property of  $\tilde\alpha(\mathbf{k})$. 

 \textbf{(II)}  \textit{For any $(i,j)\in\Sigma$, one has $\alpha((i,j))=\cT_{A_1,\theta_1}^i\circ\cT_{A_2,\theta_2}^j=\cT_{A_1^iA_2^j, ~0}=A_1^iA_2^j\times id_{\TT^s}$}.
 
  \textbf{(III)} \textit{For two linearly independent $(i, j)\in \Sigma$ and $(i',j')\in\Sigma$,  we can apply Theorem \ref{Element_Thm1} to the ergodic generators $A=A_1^iA_2^j$ and $B=A_1^{i'}A_2^{j'}$ to obtain the corresponding two positive numbers $\vep_0(A, B)$ and $\mu_0(A, B)$ for which the local rigidity holds. Let $\vep^{*}>0$ be the minimum of all such  $\vep_0(A, B)$ where  $(i, j)\in \Sigma$ and $(i',j')\in\Sigma$ are linearly independent, and let integer $\mu>0$ be the maximum of all such possible $\mu_0(A, B)$.} 

Now, we choose a sufficiently small $\vep_1$ satisfying $\vep_1<\min\left\{\frac{\delta^*}{2}, \vep^*\right\}$. Thus, when
\begin{equation*}
	d_{C^\mu}(\tilde\alpha, \alpha; M_0)<\vep_1,
\end{equation*}
and $\ZZ^2_{M_0}$ contains two linearly independent elements $\mathbf{m}, \mathbf{n}$ such that $\tilde{\alpha}(\mathbf{m})$ and $\tilde{\alpha}(\mathbf{n})$ satisfy the intersection property \ref{condIP}, we must have $\mathbf{m}, \mathbf{n}\in \Sigma$. This is due to the above property \textbf{(I)}. Write $\mathbf{m}=(m_1,m_2)$, $\mathbf{n}=(n_1,n_2)$ and set
	\[A:=A_1^{m_1}A_2^{m_2},\qquad B:=A_1^{n_1}A_2^{n_2}.\] 
By the above property \textbf{(II)},  $\alpha(\mathbf{m})=\cT_{A,0}$ and $\alpha(\mathbf{n})=\cT_{B,0}$. Note that $A$ and $B$ are still ergodic generators.
Now that $\|\tilde{\alpha}(\mathbf{m})-\TA\|_{C^\mu}<\vep_1$, $\|\tilde{\alpha}(\mathbf{n})-\TB\|_{C^\mu}<\vep_1$, and $\tilde{\alpha}(\mathbf{m})$ and $\tilde{\alpha}(\mathbf{n})$ satisfy condition \ref{condIP}, invoking Theorem \ref{Element_Thm1} we can find a $C^\infty$ near-identity conjugacy $H$ such that 
\begin{equation}\label{qdqd}
	H\circ \tilde{\alpha}(\mathbf{m})\circ H^{-1}= \TA=\alpha(\mathbf{m}),\qquad H\circ \tilde{\alpha}(\mathbf{n})\circ H^{-1}= \TB=\alpha(\mathbf{n}).
\end{equation}
This also implies that restricted on the subgroup $\mathbf{m} \ZZ+\mathbf{n}\ZZ\subset\ZZ^2$, $\tilde{\alpha}$ is conjugate to $\alpha$ via $H$.

\noindent\textbf{\large Part 2.} We still need to show that for any $\mathbf{k}\in \ZZ^2$, $\tilde{\alpha}(\mathbf{k})$ can be conjugated to $\alpha(\mathbf{k})$. In fact, we only need to verify it for the generators  $\mathbf{e}_1=(1,0)$ and $\mathbf{e}_2=(0,1)$. 

As we will see below, in general $H\circ \tilde{\alpha}(\mathbf{e}_i)\circ H^{-1}\neq \alpha(\mathbf{e}_i)$. To handle it, our plan is to construct another conjugacy which conjugates $H\circ \tilde{\alpha}(\mathbf{e}_i)\circ H^{-1}$ to $ \alpha(\mathbf{e}_i)$.

Let us write
 \[H\circ\tilde{\alpha}(\mathbf{e}_1)\circ H^{-1}=\alpha(\mathbf{e}_1)+P,\] where $P(x,y)=(P_1(x,y), P_2(x,y))$ with $P_1\in C^\infty(\TT^{d}\times\TT^s,\RR^d)$ and $P_2\in C^\infty(\TT^{d}\times\TT^s,\RR^s)$. By the commutation relation $\tilde{\alpha}(\mathbf{e}_1)\circ\tilde{\alpha}(\mathbf{m})=\tilde{\alpha}(\mathbf{m})\circ\tilde{\alpha}(\mathbf{e}_1)$ and \eqref{qdqd}, it follows that
\begin{equation*}
	P_1(Ax,y)=AP_1(x,y),\qquad P_2(Ax,y)=P_2(x,y).
\end{equation*}
As $A$ is ergodic, we obtain $P_1=0$, and $P_2(x,y)=f(y)$ is a function independent of $x$. Hence, 
\begin{equation}\label{dsfaw1}
	H\circ \tilde{\alpha}(\mathbf{e}_1)\circ H^{-1}(x,y)=(A_1x, F(y))
\end{equation}
where $F(y)=R_{\theta_1}+f(y)\in \textup{Diff}^\infty(\TT^s)$. Similarly, we can prove that  
\begin{equation}\label{dsfaw2}
	H\circ\tilde{\alpha}(\mathbf{e}_2)\circ H^{-1}(x,y)=(A_2x, G(y))
\end{equation}
where $G(y)=R_{\theta_2}+g(y)\in \textup{Diff}^\infty(\TT^s)$ with $g\in C^\infty(\TT^s,\RR^s)$.

Since $\tilde{\alpha}(\mathbf{m})=\tilde{\alpha}(m_1\mathbf{e}_1+m_2\mathbf{e}_2)=\big(\tilde{\alpha}(\mathbf{e}_1)\big)^{m_1}\circ\big(\tilde{\alpha}(\mathbf{e}_2)\big)^{m_2}$, it follows from \eqref{dsfaw1}--\eqref{dsfaw2} that
\begin{align*}
H\circ\tilde{\alpha}(\mathbf{m})\circ H^{-1}=(A_1^{m_1}A_2^{m_2}x, F^{m_1}\circ G^{m_2}(y))=(Ax, F^{m_1}\circ G^{m_2}(y)).
\end{align*}
Similarly, we also have
\begin{align*}
H\circ\tilde{\alpha}(\mathbf{n})\circ H^{-1}=(A_1^{n_1}A_2^{n_2}x, F^{n_1}\circ G^{n_2}(y))=(Bx,F^{n_1}\circ G^{n_2}(y)).
\end{align*}
Combined with \eqref{qdqd}, we get
$
	F^{m_1}\circ G^{m_2}=id_{\TT^s}$ and $F^{n_1}\circ G^{n_2}=id_{\TT^s}$. 
This also implies that
\begin{equation}\label{idTs}
	F^{q}=id_{\TT^s},\qquad G^{q}=id_{\TT^s}.
\end{equation}
with $q=|m_1n_2-m_2n_1|=|\mathbf{m}\times \mathbf{n}|$, and the integer $q>0$. Moreover, $q\theta_1\in \ZZ^s$ and $q\theta_2\in \ZZ^s$ because $m_1\theta_1+m_2\theta_2\in \ZZ^s$ and  $n_1\theta_1+n_2\theta_2\in \ZZ^s$. 

Since $d_{C^\mu}(\tilde\alpha, \alpha; M_0)<\vep_1$, by letting $\vep_1$ sufficiently small if necesary, we obtain the following Claim 1 and Claim 2.

\textbf{Claim 1:} \textit{There exists a near-identity conjugacy $V\in \textup{Diff}^\infty(\TT^s)$ such that 
\[V\circ F\circ V^{-1}=R_{\theta_1}.\]
In addition, $\|V-id_{\TT^s}\|_{C^1}\leq C \|F-R_{\theta_1}\|_{C^1}$
}

In fact, since $F$ is sufficiently $C^1$ close to $R_{\theta_1}$ and $F^q=id_{\TT^s}$ holds (see \eqref{idTs}), Claim 1 follows immediately from Lemma \ref{lem_periodid}.

\textbf{Claim 2:} \textit{We set $\tilde G:=V\circ G\circ  V^{-1}$. Then, there exists a near-identity conjugacy $\tilde{V}\in \textup{Diff}^\infty(\TT^s)$ such that 
\[\tilde V\circ \tilde G\circ  \tilde V^{-1}=R_{\theta_2},\qquad\tilde V\circ R_{\theta_1}\circ \tilde V^{-1}=R_{\theta_1}.\]}

Note that $\tilde G=R_{\theta_2}+\tilde g$ satisfies $\|\tilde g\|_{C^1}\ll 1$, and by \eqref{idTs}, $\tilde G^q=id_{\TT^s}$. Then,  invoking Lemma \ref{lem_periodid} we can find a near-identity conjugacy $\tilde V(y)=id_{\TT^s}+\tilde v(y)$ with $\tilde v\in C^\infty(\TT^s,\RR^s)$ of the form
\begin{equation*}
	\tilde v=\frac{1}{q}\sum_{i=0}^{q-2}(q-i-1)\,\tilde g\circ \tilde{G}^{i},
\end{equation*} 
such that $\tilde V\circ \tilde G\circ  \tilde V^{-1}=R_{\theta_2}$. 

Since $F\circ G=G\circ F$, it yields $R_{\theta_1}\circ\tilde G=\tilde G\circ R_{\theta_1}$, and hence $\tilde g=\tilde g\circ R_{\theta_1}$. This leads to 
\begin{align*}
	\tilde v\circ R_{\theta_1}=\frac{1}{q}\sum_{i=0}^{q-2}(q-i-1)\,\tilde g\circ \tilde{G}^{i}\circ R_{\theta_1}=\frac{1}{q}\sum_{i=0}^{q-2}(q-i-1)\,\tilde g\circ R_{\theta_1}\circ \tilde{G}^{i}
	=&\frac{1}{q}\sum_{i=0}^{q-2}(q-i-1)\,\tilde g\circ \tilde{G}^{i}\\
	=& \tilde v,
\end{align*}
from which we easily obtain $\tilde V\circ R_{\theta_1}=R_{\theta_1}\circ\tilde V $. This finally proves Claim 2.

Now, we define a smooth conjugacy $\mathcal{V}\in \textup{Diff}^\infty(\TT^d\times\TT^s)$ as follows 
\begin{equation*}
	\mathcal{V}(x,y)=(x, \tilde V\circ V(y)).
\end{equation*} 
Using \eqref{dsfaw1}--\eqref{dsfaw2} together with the above Claim 1 and Claim 2, we obtain that
\begin{align*}
	\mathcal{V}\circ H\circ \tilde\alpha(\mathbf{e}_1)\circ  H^{-1}\circ \mathcal{V}^{-1}=&A_1\times R_{\theta_1}=\alpha(\mathbf{e}_1),\\
	\mathcal{V}\circ H\circ \tilde\alpha(\mathbf{e}_2)\circ  H^{-1}\circ \mathcal{V}^{-1}=&A_2\times R_{\theta_2}=\alpha(\mathbf{e}_2).
\end{align*}
Therefore, the action $\tilde\alpha$ is $C^\infty$-conjugate to $\alpha$ via the conjugacy $U=\mathcal{V}\circ H$ provided that $d_{C^\mu}(\tilde\alpha, \alpha; M_0)<\vep_1$. Finally, it is easy to see that there exists $\vep>0$ such that whenever $d_{C^\mu}(\tilde\alpha, \alpha)<\vep$, one has  $d_{C^\mu}(\tilde\alpha, \alpha; M_0)<\vep_1$. This completes the proof.

\end{proof}

\subsection{Proof of Theorem \ref{MainThm_1}}

For the perturbed action $\tilde\alpha=\langle \cF_1, \cF_2 \rangle$, according to our assumption for each $l=1,2$, the composition $\cF_l^{q_l}=\cF_l\circ \cdots\circ\cF_l$ possesses an invariant $d$-dimensional torus homotopic to $\TT^d\times\{0\}$ $\subset$ $\TT^d\times\TT^1$. This, combined with the volume preserving condition, implies that both $\cF_1^{q_1}$ and $\cF_2^{q_2}$ satisfy the intersection property since the fiber is of dimension one. Therefore, based on Theorem \ref{Element_Thm1} and Lemma \ref{lem_periodid}, the rest of the proof is just analogous to that of Theorem \ref{MainThm_0}.

\subsection{Proof of Theorem \ref{MainThm_2}}
Since the $\ZZ^k$ action $\rho_0: \ZZ^k\to \textup{Aut}(\TT^d)$ is higher rank on the base $\TT^d$, one can find a subgroup $\Sigma\subset \ZZ^k$ with $\Sigma\cong \ZZ^2$ such that $\rho_0(\mathbf{n})$ is ergodic on $\TT^d$ for all $\mathbf{n}\in\Sigma\setminus\{0\}$. 

For the perturbed action $\tilde\rho$, by assumption the generators  $\tilde\rho(\mathbf{e}_1),\cdots, \tilde\rho(\mathbf{e}_k)$ have a common invariant $d$-dimensional torus homotopic to $\TT^d\times\{0\}$ $\subset$ $\TT^d\times\TT^1$. This, combined with the volume preserving condition and the fiber is of dimension 1, implies that for every $\mathbf{m}\in\ZZ^k$, the map $\tilde\rho(\mathbf{m}): \TT^d\times\TT^1\to \TT^d\times\TT^1$  satisfies the intersection property. Then, by applying Theorem \ref{Element_Thm1} to the subgroup $\Sigma$, we obtain that the restriction  $\tilde\rho\big|_{\Sigma}$ of $\tilde\rho$ to the subgroup $\Sigma$ can be $C^\infty$-conjugated to $\rho\big|_{\Sigma}=\rho_0\times id_{\TT^1}\big|_{\Sigma}$ via a conjugacy $H$, provided that $\textup{dist}_{C^\mu}(\tilde\rho,\rho)$ is sufficiently small.

To complete the proof, it remains to show that for other $\mathbf{m}\in \ZZ^k\setminus\Sigma$, $\tilde\rho(\mathbf{m})$ is $C^\infty$-conjugate to $\rho(\mathbf{m})$. Indeed, we only need to verify it for the generators $\mathbf{e}_1,\cdots, \mathbf{e}_k$ of $\ZZ^k$. Here, we only check it for $\mathbf{e}_1$, and other cases are similar.

Using arguments analogous to the proof of Theorem \ref{MainThm_0}, the ergodicity of $\rho_0\big|_{\Sigma}$ on the base $\TT^d$ and the commutativity imply that
\[H\circ \tilde\rho(\mathbf{e}_1)\circ H^{-1}=\rho(\mathbf{e}_1)+(0, f(y) ) \]
for some  $f(y)\in C^\infty(\TT^1,\RR^1)$. As a result of the intersection property of $\tilde\rho(\mathbf{e}_1)$,  $f(y)$ has to be zero. Therefore, $H\circ \tilde\rho(\mathbf{e}_1)\circ H^{-1}=\rho(\mathbf{e}_1).$ The case of $\mathbf{e}_2, \cdots, \mathbf{e}_k$ can be proved in the same fashion as that of $\mathbf{e}_1$. This finally proves Theorem \ref{MainThm_2}.


\end{document}